\documentclass[12pt]{amsart}
\usepackage{amsthm, amsmath}
\usepackage{amssymb} 
\usepackage{graphicx}

\theoremstyle{plain}
\newtheorem{theorem}{Theorem}[section]
\newtheorem{lemma}[theorem]{Lemma}
\newtheorem{cor}[theorem]{Corollary}
\newtheorem{prop}[theorem]{Proposition}

\newtheorem{thr}{Theorem}
\newtheorem{propo}[thr]{Proposition}

\theoremstyle{remark}

\theoremstyle{definition}

\newtheorem{defn}[theorem]{Definition}
\newtheorem{Remark}[theorem]{Remark}
\newcommand{\NN}{\mathbb N}

\newcommand{\ZZ}{\mathbb Z}
\newcommand{\QQ}{\mathbb Q}
\newcommand{\RR}{\mathbb R}

\newcommand{\FF}{\mathbb F}

\newcommand\Adm{{\rm Adm}}
\newcommand\Perm{{\rm Perm}}
\newcommand\calA{{\mathcal A}}

\title{Alcoves associated to special fibers 
of local models}
\author{Thomas J. Haines \,\, and  \,\, Ng\^{o} Bao Ch\^au }

\date{}

\begin{document}
\maketitle

%\begin{abstract}

%\end{abstract} 

%%% Section 1
\markboth{T. Haines and B.C. Ng\^o}
{Alcoves associated to special fibers of local models}

\section{Introduction}

Let $G$ be a classical group over the $p$-adic field $\QQ_p$, and let $\mu$ be a dominant 
minuscule coweight of $G$.  We assume that the pair $(G,\mu)$ comes from a PEL-type 
Shimura datum 
$({\bf G},X,K)$ where the $p$-component $K_p$ of the compact open subgroup $K$ is an 
Iwahori subgroup of ${\bf G}(\QQ_p) = G(\QQ_p)$.  We assume that the corresponding Shimura 
variety ${\rm Sh}({\bf G},X,K)$ has a model over a $p$-adic number ring ${\mathcal O}_{\mathfrak p}$.  
In this case ${\rm Sh}({\bf G},X,K)$ is said to have {\em Iwahori-type bad reduction} at the prime 
${\mathfrak p}$.  A fundamental problem is to understand the geometry of the special fiber.  In the 
prototypical example of the modular curve $Y_0(p)$, the special fiber is a union of two smooth 
curves which intersect transversally.  In the general case the global geometry cannot be easily 
described and the singularities can be quite complicated (see \cite{RapAA}).   

The local geometry can be approached using the Rapoport-Zink local model $M_\mu$ 
\cite{Rap-Zink}.  This projective ${\mathcal O}_{\mathfrak p}$-scheme is a local model for the 
singularities in the special fiber of ${\rm Sh}({\bf G},X,K)$, but it is defined in terms of the 
pair $(G,\mu)$ using linear algebra and is somewhat easier to deal with than the Shimura 
variety itself.   

In this paper we address certain combinatorial questions which arise in the study of the 
special fiber $M_{\mu, \overline{\FF}_p}$ of $M_\mu$.  To fix ideas, consider the case where 
$G$ is {\em split}.  It turns out that $M_{\mu, \overline{\FF}_p}$ can be regarded as a 
finite-dimensional union of Iwahori-orbits in the affine flag variety for $G(\overline{\FF}_p(\!(t)\!))$,
%%(since $G$ is split we can define it over $\ZZ$, so this has a meaning)
see \cite{Goertz1} and \cite{H-N}.  It follows 
that $M_{\mu,\overline{\FF}_p}$ has a stratification indexed by a finite set ${\rm Perm}(\mu)$ of 
the extended affine Weyl group $\widetilde{W}(G)$ for $G$.  Equivalently, ${\rm Perm}(\mu)$ can 
be regarded as a finite set of alcoves in the affine Coxeter complex determined by a splitting for 
$G$.  Let $W_0$ denote the finite Weyl group of $G$.  For each translation $\lambda$ in the 
$W_0$-orbit of $\mu$ there is an Iwahori-orbit in $M_{\mu, \overline{\FF}_p}$ corresponding 
to $\lambda$; in other words, the element $t_\lambda$ in $\widetilde{W}(G)$ is contained in 
${\rm Perm}(\mu)$.  We let ${\rm Adm}(\mu)$ denote the subset of ${\rm Perm}(\mu)$ indexing 
those strata which lie in the closure of the stratum indexed by $t_\lambda$, for some 
$\lambda \in W_0(\mu)$.  

We are primarily concerned with {\em the equality ${\rm Adm}(\mu) = {\rm Perm}(\mu)$}, which has 
important geometric content.  For example, if the ${\mathcal O}_{\mathfrak p}$-scheme $M_\mu$ 
is {\em flat}, this equality automatically holds.  On the other hand, the equality has 
been established combinatorially when $G = {\rm GL}_n$ or ${\rm GSp}_{2n}$ by Kottwitz 
and Rapoport \cite{KoRa}, and this was exploited by U. G\"{o}rtz in his proof of 
the flatness of the local models associated to $\mbox{Res}_{E/F}{\rm GL}_n$, where $E/F$ 
is an unramified extension of $p$-adic fields. 

In \cite{KoRa}, Kottwitz and Rapoport propose purely combinatorial 
definitions of the sets ${\rm Adm}(\mu)$ and 
${\rm Perm}(\mu)$, for every split group $G$ and every dominant 
coweight $\mu$.  An element $x \in \widetilde{W}(G)$ belongs to ${\rm Adm}(\mu)$ if 
$x \leq t_\lambda$ in the Bruhat order on $\widetilde{W}(G)$, for some $\lambda \in W_0(\mu)$.  
The element $x$ belongs to ${\rm Perm}(\mu)$ if the vertices of the corresponding alcove 
satisfy certain inequalities similar to those arising from the definition of local models; see 
section 3 for details.  (These definitions agree with the geometric interpretations above when 
$G = {\rm GL}_n$ or ${\rm GSp}_{2n}$ and $\mu$ is minuscule, see loc. cit..)  
Kottwitz and Rapoport show that {\em the inclusion ${\rm Adm}(\mu) \subset {\rm Perm}(\mu)$ 
is always valid.} We adopt their point of view and assume from now on that $G$ is split.  
We address the question of whether the equality ${\rm Adm}(\mu) = {\rm Perm}(\mu)$ holds 
in full generality, in particular for {\em non-minuscule} coweights $\mu$.  

The main result of this paper is that the situation is as nice as possible for groups of type 
$A_{n-1}$.

\begin{thr} For any root system of type $A_{n-1}$, the equality 
$\Adm(\mu) = \Perm(\mu)$ holds for every dominant coweight $\mu$. 
\end{thr}

Actually, we derive this theorem from a more general result available for 
any root system. We introduce the new notion of the set of 
$\mu$-{\em strongly permissible} alcoves $\Perm^{st}(\mu)$.  
Like the Kottwitz-Rapoport notion of $\Perm(\mu)$, this is a set 
of alcoves determined by imposing conditions vertex-by-vertex. 
Although the condition we require on vertices might look more technical than
that of Kottwitz-Rapoport, a similar notion for finite Weyl groups 
already occurred in a classical theorem of Bernstein-Gelfand-Gelfand
\cite{BGG1}. We prove the following statement.

\begin{propo}
For any root system  $R$, the inclusion $\Adm(\mu)\supset \Perm^{st}(\mu)$ 
holds for every dominant coweight $\mu$.
\end{propo}

When the root system $R$ is of type $A_{n-1}$, we prove that 
$\mu$-permissibility and $\mu$-strong permissibility are equivalent and
derive Theorem 1 from Proposition 2 and the inverse inclusion proved 
by Kottwitz-Rapoport. 

In this paper we introduce some notions of cones, acute or obtuse,
in the set of alcoves. These notions play an important role in the proof 
of Proposition 2 and may also be useful for other combinatorial 
questions.  Another ingredient is a lemma due to Deodhar \cite{Deod1} 
which is also used in \cite{Deod2}.

The equality $\Adm(\mu)=\Perm(\mu)$ turns out to be false in general. 
 
\begin{thr}
If $R$ is an irreducible root system of rank $\,\, \geq 4$ and not of type $A_{n-1}$,
then $\Adm(\mu)\not=
\Perm(\mu)$ for every sufficiently regular dominant coweight $\mu$.
\end{thr}

Theorem 3 is based on counter-examples due to Deodhar \cite{Deod2}.  
Deodhar determines which finite Weyl groups 
$W$ have the property that $w \le w'$ in the Bruhat order if and only if 
$w(\lambda) - w'(\lambda)$ is a sum of positive 
coroots for every dominant coweight $\lambda$.  
He proves that $W$ 
has this property if and only if the irreducible components of the 
associated root system $R$ are of type $A_{n-1}$ or of rank $\,\, \leq 3$.   
Our proof of Theorem 3 uses the counter-examples he
gives explicitly in the other cases.
However Theorem 1 is proved independently and in fact in section 8 we give 
an alternate, perhaps more conceptual, proof of Deodhar's theorem 
for root systems of type $A_{n-1}$.

In light of Theorem 3, and motivated by the study of local models attached to 
certain {\em nonsplit} groups, Rapoport proposed that the equality ${\rm Adm}(\mu) 
= {\rm Perm}(\mu)$ might be valid for all root systems as long as $\mu$ is a {\em sum 
of minuscule coweights}.  In section 10 we show this is indeed the case for the symplectic group.

\begin{thr}
Let $\mu$ be a sum of minuscule coweights for the group ${\rm GSp}_{2n}$.  
Then the equality ${\rm Adm}(\mu) = {\rm Perm}(\mu)$ holds.
\end{thr}

Theorem 4 relies on a description of the root system of ${\rm GSp}_{2n}$ as the 
``fixed-point'' root system $R^{[\Theta]}$ with respect to the nontrivial automorphism 
$\Theta$ of the root system $R$ for ${\rm GL}_{2n}$, \`{a} la Steinberg \cite{Steinberg}.  
This idea was also exploited in \cite{KoRa}.  In fact in section 9 we prove the following statement is valid for every 
dominant coweight, at least in the case of the symplectic group.
%%In fact the general results proved about admissible sets and automorphisms in section 9 yield,  at least in the case of the symplectic group, the following criterion for the equality of the admissible and permissible sets.  Theorem 4 is a consequence of this more precise result.    

\begin{propo}    
 The equality ${\rm Adm}^\Theta(\mu) ={\rm Perm}(\mu) \cap \widetilde{W}({\rm GSp}_{2n})$ holds
for any dominant coweight $\mu$ of the root system $R^{[\Theta]}$ for ${\rm GSp}_{2n}$.
\end{propo}

Here the sets ${\rm Adm}^\Theta(\mu)$ and ${\rm Perm}^\Theta(\mu)$ are the subsets of 
$\widetilde{W}({\rm GSp}_{2n})$ analogous to the subsets ${\rm Adm}(\mu)$ and 
${\rm Perm}(\mu)$ of $\widetilde{W}({\rm GL}_{2n})$.  
%%There is a similar criterion for the equality of admissible and strongly permissible sets 
Moreover, the equality ${\rm Perm}^\Theta(\mu)={\rm Perm}(\mu) \cap \widetilde{W}({\rm GSp}_{2n})$
can be proved in the case where $\mu$ is a sum of minuscule coweights by a method of 
Kottwitz-Rapoport \cite{KoRa}, yielding Theorem 4. It is also worthwhile to remark that the set ${\rm Perm}(\mu) \cap \widetilde{W}
({\rm GSp}_{2n})$ has itself some geometric meaning for every $\mu$. This set parametrizes 
the strata in the special fiber of the local model considered in \cite{H-N}, since the definition 
thereof uses the standard representation of the symplectic group.

Theorem 1 is expected to play a role in proving the flatness of the 
Pappas-Rapoport local models attached to the group $\mbox{Res}_{E/F}{\rm GL}_n$, 
where $E/F$ is a totally ramified extension of $p$-adic fields.  For details we 
refer the reader to forthcoming work of U. G\"{o}rtz \cite{Goertz2}. Likewise, 
Theorem 4 (or Proposition 5) is expected to play a role in proving the flatness of analogous local 
models attached to $\mbox{Res}_{E/F}{\rm GSp}_{2n}$.  
%%(Indeed, a minuscule cocharacter for the latter can be identified with a sum of minuscule coweights for ${\rm GSp}_{2n}$; similar remarks apply to $\mbox{Res}_{E/F}{\rm GL}_n$.)
Thus, our study of ${\rm Adm}(\mu)$ and ${\rm Perm}(\mu)$ for {\em non-minuscule} 
coweights $\mu$ of split groups has ramifications for Shimura varieties attached to 
certain {\em non-split} groups.    
    
For applications to the bad reduction of PEL-type Shimura varieties attached to orthogonal groups 
one would like to determine whether the analogs of Theorem 4 and Proposition 5 hold for the split orthogonal 
groups ${\rm O}(2n)$ (PEL Shimura varieties arise for groups of type $A$,$\,\, C$, and $D$).  The methods of this paper do not seem to give much information in that situation.  However, in section 11 we show that the analog of Proposition 5 for the odd orthogonal groups does not hold.  This can be understood in terms of a {\em non-inheritance} property of the Bruhat order: the group ${\rm SO}(2n+1)$ can be realized as the fixed point group ${\rm SL}(2n+1)^\Theta$ for a certain involution $\Theta$, and this gives a corresponding embedding of affine Weyl groups $W_{\rm aff}(B_n) \hookrightarrow W_{\rm aff}(A_{2n})$.  However, the Bruhat order on the former is not inherited from the Bruhat order on the latter (in contrast to the symplectic case, cf. Proposition 9.6).

\section{Affine Weyl group attached to a root system}

Let us fix notation and recall basic facts about affine Weyl group and 
alcoves. For proofs, we refer to Humphreys' book \cite{Hum}, especially chapter 4. 

Let $(X^*,X_*,R,\Check{R})$ be a root system. We assume throughout 
this paper that the root system is {\em reduced and irreducible}.  When there is no chance 
of confusion we will denote the root system $(X^*, X_*, R, \Check{R})$ 
simply by $R$. Let $\Pi$ be a base of $R$.  Let $R^+$ (resp. $R^-$) 
denote the set of positive (resp. negative) roots. The cardinality of $\Pi$ 
will be called the {\em rank} of $R$.  Denote by $\langle ~ , ~ \rangle : X^* \times X_* \rightarrow \ZZ$ the perfect pairing making $X^*$ and $X_*$ dual free abelian groups.

Corresponding to $\alpha \in R$ we have a reflection
$s_{\alpha}$, acting on $V=X_*\otimes\RR$ by 
$s_{\alpha}(x) = x - \langle \alpha , x \rangle\check{\alpha}$.
The Weyl group $W_0$ is the subgroup of ${\rm GL}(V)$ 
generated by these reflections. It is known that $W_0$ is generated by
$S = \lbrace s_{\alpha} ~ | ~ \alpha 
\in \Pi \rbrace$ as a finite Coxeter group.  

Let $H_\alpha$ denote the hyperplane ({\em wall}) fixed by the reflection
$s_\alpha$, for every $\alpha \in R$.  The connected components of the set
$V - \bigcup_{\alpha \in R} H_\alpha$ will be called {\em chambers}.  
The finite Weyl group $W_0$ acts simply transitively on the set $\mathcal C$
of chambers. There is a distinguished chamber 
$$C_0 = \{ x \in V \,\, | \,\, \langle \alpha , x \rangle > 0, \,\,\, \mbox{for every}  
\,\,\, \alpha \in \Pi \}$$
which will be called the {\em dominant chamber}.

Corresponding to $\alpha \in R$ and $k\in\ZZ$ we have an {\em affine} reflection
$s_{\alpha,k}$, acting on $V=X_*\otimes\RR$ by 
$s_{\alpha,k}(x) = x - (\langle \alpha , x \rangle-k)\check{\alpha}$.
The affine Weyl group $W_{\rm aff}$ is the subgroup of ${\rm Aut}(V)$ 
generated by these reflections. It is known that $W_{\rm aff}$ is generated,
as a Coxeter group, by $S_{\rm aff}=S\cup\lbrace s_{\tilde{\alpha},1} \rbrace$
where $\tilde{\alpha}$ is the unique highest root of $R$. 
Moreover, $W_{\rm aff}$ is the semi-direct product $W_{\rm aff}=
W_0\ltimes \Check Q$ where $\Check Q$ is the lattice generated by 
$\Check R$ acting on $V$ by translation.  

Let us denote $\tilde R=R\times\ZZ$.  Let $H_{\alpha,k}$ denote the hyperplane in
$V$ fixed by the reflection $s_{\alpha,k}$ for every $(\alpha,k)\in\tilde R$.  
The connected components of the set 
$V - \bigcup_{(\alpha,k) \in \tilde R} H_{\alpha,k}$ will be called {\em alcoves}.
The affine Weyl group $W_{\rm aff}$ acts simply transitively on the set $\mathcal A$ 
of the alcoves. There is a distinguished alcove 
$$A_0 = \{ x \in C_0 \,\, | \,\, \langle \tilde{\alpha} , x \rangle < 1 \}$$
that will be called the {\em base alcove}. 

The extended Weyl group $\widetilde W=W_0\ltimes X_*$ also acts on the set  
$\mathcal A$. Indeed for any coweight $\lambda\in X_*$ the {\em translation} 
$t_\lambda$ by the vector $\lambda$ sends a wall on another wall since
$\langle\alpha,\lambda\rangle\in\ZZ$ for every $\alpha\in R$, therefore 
sends an alcove to another alcove. 
Let $\Omega$ be the isotropy group in $\widetilde W$ of the base alcove. It is known 
that $\widetilde W=W_{\rm aff}\rtimes \Omega$. With the help of this decomposition, the 
Bruhat order and the length function on the Coxeter group $W_{\rm aff}$  can be 
extended to $\widetilde W$, which is not a Coxeter group in general. 
For $w,w'\in W_{\rm aff}$,
$\tau,\tau'\in\Omega$, we say $w\tau \leq w'\tau'$ if and only if $w\leq w'$
and $\tau=\tau'$. We put $l(w\tau)=l(w)$. 

Let us recall some basic facts on minimal galleries.
A gallery of length $l$ is a sequence of alcoves $A'=A'_0,\ldots,A'_l=A''$ such 
that $A'_{i-1}$ and $A'_{i}$ share a wall $H_i$, for $i = 1,\dots,l$. 
The gallery is minimal if there 
does not exist a gallery going from $A'$ to $A''$ with length strictly less than $l$.
In general there exist more than one {\em minimal gallery} going from $A'$ from 
$A''$, but the set of walls $H_i$ depends only on $A'$ and $A''$: it is the set of 
walls $H=H_{\alpha,k}$ separating $A'$ and $A''$, i.e., those such that 
$A'$ and $A''$ lie in different connected components of $V-H$.

Minimal galleries are very closely related to reduced expressions in the Coxeter
system $(W_{\rm aff},S_{\rm aff})$.
Let $x\in W_{\rm aff}$ and $x=s_1s_2\ldots s_l$ be a {\em reduced expression} with 
$s_1,\ldots,s_l\in S_{\rm aff}$,  i.e., such an expression with minimal length $l=l(x)$.
Let $A_0,A_1,\ldots,A_l$ be the alcoves defined by  
$A_{i}=s_1s_2\ldots s_i A_0$. Obviously, $A_{i-1}$ and $A_{i}$ share a wall $H_i$.
Hence, this sequence of alcoves  forms a gallery going from the base alcove 
$A_0$ to $A_l=x A_0$. In this way, the reduced expressions
of $x$ correspond bijectively to the minimal galleries going from $A_0$ to $xA_0$.
The cardinality of the set of hyperplanes $H_{\alpha,k}$ with
$(\alpha,k)\in\tilde R$  separating $A_0$ and $xA_0$ is equal to the length
$l(x)$.

Another basic fact that will be used in the sequel is the following. If an expression 
$x=s_1s_2\ldots s_l$ is not reduced, there are integers $i<j$ such that we can delete 
$s_i$ and $s_j$ without changing $x$:
$$x=s_1\ldots\hat s_i\ldots \hat s_j \ldots s_l.$$
Equivalently, a gallery $A'_0,\ldots,A'_l$ is not minimal 
if there exist integers $i<j$ such that $H_i=H_j$, where for every $i$, $H_i$ denotes the wall shared by the consecutive alcoves $A'_{i-1}$ and $A'_i$.

\section{The $\mu$-admissible and $\mu$-permissible sets}

In the sequel, a dominant coweight $\mu$ will be fixed 
and we will denote by $\tau$ the unique element of $\Omega$ 
such that $t_\mu\in W_{\rm aff}\tau$.
 Let ${\rm Conv}(\mu)$ denote the convex hull of the set 
$\lbrace \lambda ~ | ~ \lambda \in W_0(\mu) \rbrace$.  
In \cite{KoRa}, R. Kottwitz and M. Rapoport introduce the notions of the 
$\mu$-{admissible} and $\mu$-{permissible} subsets 
of $\widetilde W$. 

\begin{defn} ({\bf Kottwitz-Rapoport})
\begin{enumerate}
\item Let $\Adm(\mu)=\lbrace x \in \widetilde{W}\mid
x \leq t_\lambda \,\, 
\mbox{for some} \,\, \lambda \in W_0(\mu) \rbrace$.  An element $x \in \Adm(\mu)$
 is called $\mu$-{\em admissible}.

\item Let $\Perm(\mu)$ be the set of elements $x \in W_{\rm aff}\tau$
such that $x(a) - a \in {\rm Conv}(\mu)$ for every vertex $a\in A_0 $. 
An element $x \in \Perm(\mu)$ is called $\mu$-{\em permissible}. 
\end{enumerate}
\end{defn}

We remark that in general the facets of $A_0$ having minimal dimension are not points, so our use of the word ``vertex'' in the definition above may be considered inaccurate.  Throughout this paper we shall ignore this subtlety and continue to imagine that the facets of $A_0$ of minimal dimension are in fact points.  It is not difficult to translate our resulting arguments and statements into ones which do not abuse 
terminology.

The following theorem is proved in \cite{KoRa}.  

\begin{theorem} ({\bf Kottwitz-Rapoport})

\begin{enumerate}
\item
For any root system $R$ and dominant coweight $\mu$, 
$\Adm(\mu) \subset \Perm(\mu)$.
\item
Let $R$ be the root system attached to ${\rm GL}_n$ or ${\rm GSp}_{2n}$.  
Let $\mu$ be a minuscule coweight of $R$.  Then $\Adm(\mu) = \Perm(\mu)$.
\end{enumerate}
\end{theorem}

Recall that a coweight $\mu$ is {\em minuscule} if 
$\langle \alpha, \mu \rangle \in \lbrace -1,0,1 \rbrace$, for every root 
$\alpha$.  For root systems of type $A_{n-1}$, the second statement can be 
generalized to every 
dominant coweight $\mu$.

\begin{theorem} For any root system of type $A_{n-1}$, the equality 
$\Adm(\mu) = \Perm(\mu)$ holds for every dominant coweight $\mu$. 
\end{theorem}
 
According to Kottwitz-Rapoport, it is sufficient to prove the inclusion
$\Adm(\mu)\supset \Perm(\mu)$. 

\section{Obtuse cones and $\mu$-strongly permissible sets}

Let $B_0$ denote the negative obtuse cone in $V$ generated by the coroots 
$-\Check\alpha$ with $\alpha\in R^+$.
The convex hull ${\rm Conv}(\mu)$ occurring in the definition of the 
$\mu$-permissible set may be usefully described as the intersection of $|W_0|$ obtuse cones
$${\rm Conv}(\mu)=\bigcap_{w\in W_0} w\mu+ w(B_0).$$
We can rephrase the definition of $\mu$-permissibility by saying that
 $x\in W_{\rm aff}\tau$ is $\mu$-permissible if and only if  
for every vertex $v$ of the base alcove and for any $w\in W_0$, we have
$$x(v)\in t_{w\mu}(v)+w(B_0).$$
The notion of $\mu$-strong permissibility consists in requiring a little bit more
than the last inclusion.

Let $v$ be an element of $V$, $W_{\rm aff}(v)$ the orbit of $v$ under the action 
of the affine Weyl group $W_{\rm aff}$. In the spirit of a paper of 
Bernstein-Gelfand-Gelfand \cite{BGG1}, it seems reasonable to consider the subset $B(v,w)$
of $W_{\rm aff}(v)$ defined as follows.

\begin{defn}
For every $v\in V$ and $w\in W_0$, let $B(v,w)$ be the set of elements
$v'$ of the form $s_r \ldots s_1(v)$
with $s_i=s_{\alpha_i,k_i}$ with $(\alpha_i,k_i)\in\tilde R$ such that for 
all $i=1,\ldots,r$, 
$$s_i s_{i-1}\ldots s_1(v)\in s_{i-1}\ldots s_1(v)+ w(B_0).$$ 
\end{defn}

We obviously have the inclusion 
$$
B(v,w)\subset W_{\rm aff}(v)\cap (v+w(B_0)).
$$  
In contrast to the notion of $\mu$-permissibility, we use these smaller sets $B(v,w)$ to define the notion of $\mu$-strong permissibility.

\begin{defn}
For every dominant coweight $\mu\in X_*$, let ${\Perm}^{st}(\mu)$ be the set 
of $x\in W_{\rm aff}\tau$ such that for every vertex $a$ of $A_0$ and every $w\in W_0$
the element $x(a)$ lies in $B(t_{w\mu}(a), w)$. Such an element $x$ is said to be
$\mu$-{\em strongly permissible}.
\end{defn}
We obviously have the inclusion ${\Perm}^{st}(\mu) \subset \Perm(\mu)$.
We will prove the following stronger statement.

\begin{prop}
For any root system $R$ and dominant coweight $\mu\in X_*$, the inclusion
$\Adm(\mu)\supset \Perm^{st}(\mu)$ holds.
\end{prop}

To prove this theorem we need the notion of acute cone to be discussed
in section 5. The proof itself is postponed to section 6.

The following lemma shows that in the case of a root system of type $A_{n-1}$
the notions of $\mu$-permissibility and $\mu$-strong permissibility coincide.
Thus Theorem 3.3 follows from Proposition 4.3.

\begin{lemma}
Suppose $R$ is of type $A_{n-1}$. Let $v$ be a vertex of some alcove $A$.
For every $w\in W_0$, the equality 
$$B(v,w)=W_{\rm aff}(v)\cap (v+w(B_0)) $$
holds. In particular $\Perm(\mu)=\Perm^{st}(\mu)$ for every dominant 
coweight $\mu$. 
\end{lemma}

\begin{proof}
For a root system $R$ of type $A_{n-1}$, the vertices of alcoves belong
to the coweight lattice $\Check P=\lbrace v\in V\mid \langle \alpha, v \rangle\in\ZZ\ 
\mbox{for every root}\ \alpha\in R \rbrace$. Moreover, if $v\in \Check P$, the
orbit $W_{\rm aff}(v)$ is exactly the set of elements $v'\in\Check P$ such that 
$v'-v\in \Check Q$, as one sees using the relation $t_{\check{\alpha}} = s_{\alpha, 1}s_\alpha$, for any $\alpha \in R$.
Since $\Check Q$ admits $\lbrace w(\Check\alpha) \mid \alpha\in\Pi \rbrace$
as $\ZZ$-basis, any element
$$v'\in W_{\rm aff}(v)\cap (v+w(B_0))$$
can be written uniquely in the form
$v'=v-\sum_{\alpha\in\Pi} n_\alpha w(\Check\alpha)$
with $n_\alpha\in\NN$. Thus, it is sufficient to prove the inclusion 
$v-w(\Check\alpha)\in B(v,w)$ for any simple root $\alpha\in\Pi$.
Let $k$ be the integer $\langle w(\alpha),v\rangle$.  Obviously, 
$v-w(\Check\alpha)=s_{w(\alpha),k-1}(v)$ and we are done.
 \end{proof}

\includegraphics{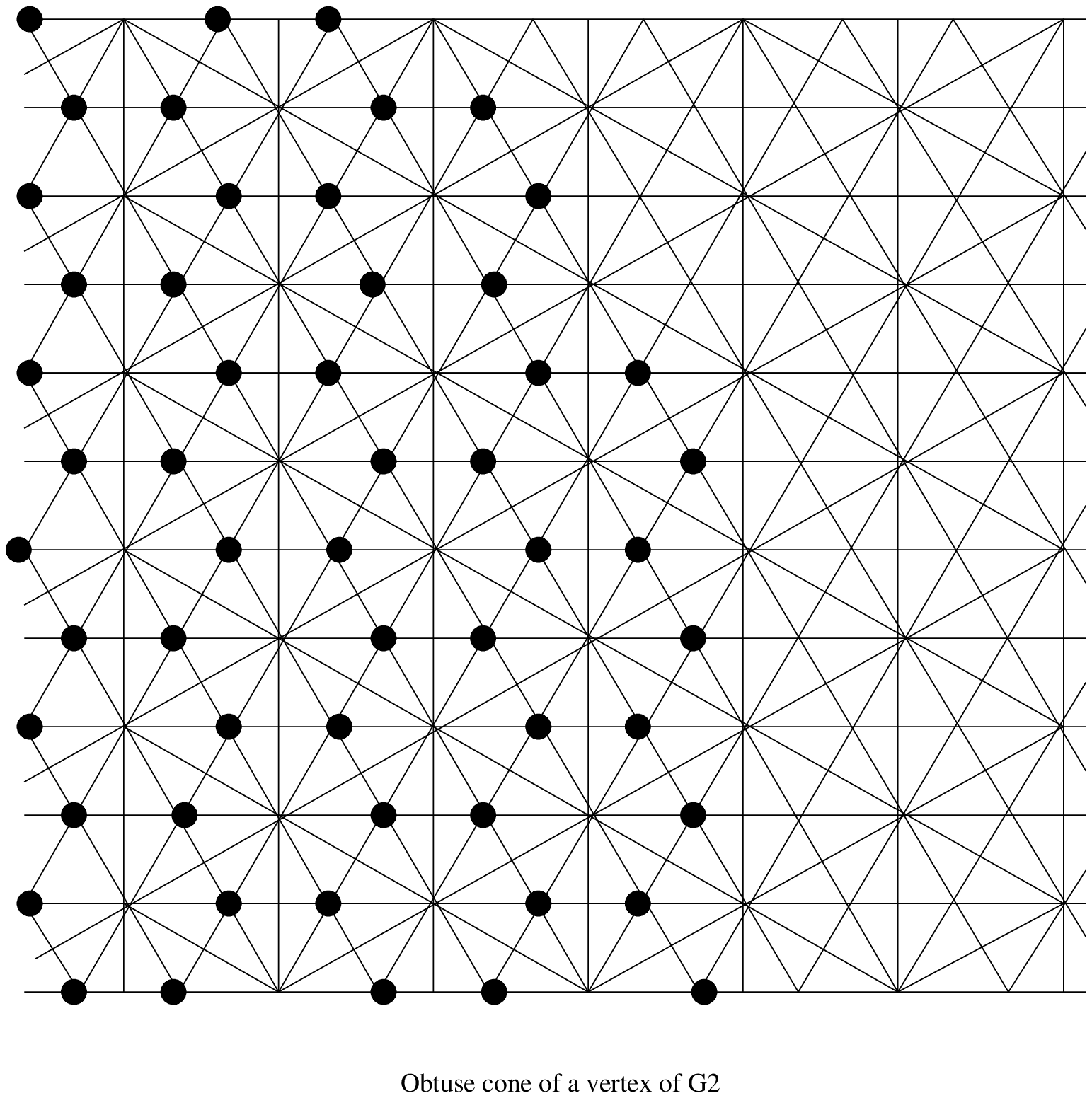}

\begin{Remark}
Let $R$ be any root system and let $v$ be a {\em special} vertex for $R$, i.e. an element of $\check{P}$.  Then for every element 
$w \in W_0$ we have inclusions
$$
\lbrace v - \sum_{\alpha \in \Pi} \NN w(\check{\alpha}) \rbrace ~ \subset ~ B(v,w) ~ 
\subset ~ W_{\rm aff}(v) \cap (v+ w(B_0)) ~ \subset ~ \lbrace v - \sum_{\alpha \in \Pi} \NN w(\check{\alpha}) \rbrace,
$$
the first and last inclusions being proved as in the proof of Lemma 4.4; 
in particular the displayed equality in Lemma 4.4 holds 
provided that $v$ is a special vertex.  
The importance of the first set above stems from its interpretation as the 
coweights appearing in the character of the Verma module 
attached to the vertex $v$ and the system of positive roots $w(R^+)$.  
One can see from the example of nonspecial vertices for the
root systems $B_2$ and $G_2$ that in general $\lbrace v - \sum_{\alpha \in \Pi} \NN w(\check{\alpha}) \rbrace \nsupseteq B(v,w)$.  
One might hope that the second inclusion 
is always an equality (and thus that $\Perm^{st}(\mu)$ always equals 
$\Perm(\mu)$).   However, this also fails in general, as we shall see in 
$\S 7$ using the examples of Deodhar \cite{Deod2}.  Of the three sets above, $B(v,w)$ seems to be the 
most related to the Bruhat order but also the 
most difficult to be visualized.

The equality of the sets in Lemma 4.4 for {\em every} vertex imposes 
rather strong conditions on the underlying root system $R$.  In fact, 
as we shall see in section 8 (Remark 8.3), if $w_0$ is the longest element of $W_0$ and
$$
B(v,w_0) = W_{\rm aff}(v) \cap (v + w_0(B_0))
$$ 
for every vertex $v$, then $R$ is necessarily of type $A_{n-1}$ or of rank $\leq 3$.
\end{Remark}

\section{Acute cones}

Let $\bar C_0$ be the closure of the dominant  chamber $C_0$.
For any $x\in V$ and $w\in W_0$, the acute cone
$\bar C(x,w)=x+w(\bar C_0)$ is the translation of $w(\bar C_0)$ by the vector $x$. 
For $A\in\calA$ and $w\in W$, we want to define a subset $C(A,w)$ 
of $\calA$ which is as close as possible to a cone $\bar C(x,w)$.

The easiest way to do so consists in choosing a point $a$ inside the 
base alcove $A_0$. Since $W_{\rm aff}$ acts simply transitively on $\calA$,
each alcove $A$ contains a unique point $A(a)$ conjugate to $a$.  

\begin{defn}
$C_a(A,w) = \lbrace A' \in {\mathcal A} ~ | ~ A'(a)  \in A(a)+ w(\bar C_0) \rbrace$.
\end{defn}

The disadvantage of this definition is that it really depends on the choice of 
$a\in A_0$. In order to give a more intrinsic definition, we need the notion
of direction of a gallery.

Let $w\in W_0$. A wall $H=H_{\alpha,k}$ can also given by $H_{-\alpha,-k}$.
But only one root from $\lbrace\alpha,-\alpha\rbrace$ lies in $w(R^+)$; 
let us assume $\alpha\in w(R^+)$. We define the $w$-positive half-space
$H^{w+}$ by $H^{w+}=\lbrace v\in V\mid \langle\alpha,x\rangle >k \rbrace$. 
Let $H^{w-}$ be the other half-space.

\begin{defn}
A gallery $A'_0,\ldots,A'_l$ is said to be in the {\em $w$-direction} if for any 
$i=1,\ldots,l$, letting $H_i$ denote the common wall of $A'_{i-1}$ and of $A'_{i}$,
the alcove $A'_{i-1}$ lies in the $w$-negative half-space $H_i^{w-}$ and 
the alcove $A'_{i}$ lies in the $w$-positive half-space $H_i^{w+}$.
\end{defn}

Obviously, the gallery $A'_0,\ldots,A'_l$ is in the $w$-direction if and only if the
inverse gallery $A'_l,\ldots,A'_0$ is in the $ww_0$-direction where $w_0$ is the
maximal length element of $W_0$. We will sometimes say the gallery 
$A'_l,\ldots,A'_0$ is in the {\em $w$-opposite direction}.

\includegraphics{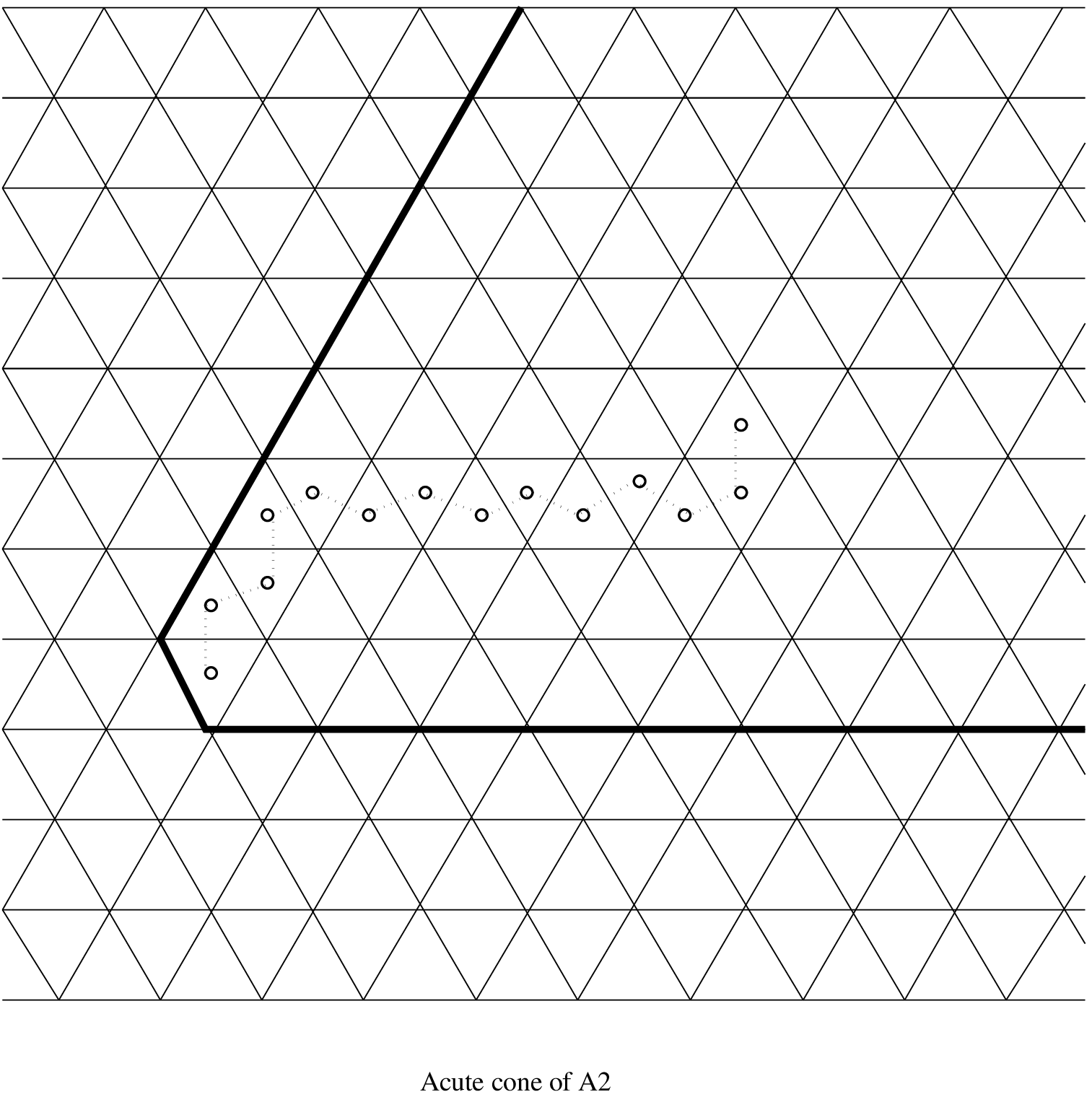}

\begin{lemma}
Any gallery $A'_0,\ldots,A'_l$ in some $w$-direction is 
automatically minimal. If there exists a gallery going from $A'_0$ 
to $A'_l$ in the $w$-direction, then any other minimal gallery going
from $A'_0$ to $A'_l$ is in the $w$-direction too.
\end{lemma}

\begin{proof}
Suppose the gallery $A'_0,\ldots,A'_l$ is not minimal.  This means there exist
two integers $i<j$ such that $H_i=H_j$. Let $H$ denote this hyperplane and 
suppose $H_k\not= H$ for any integer $k$ in between $i<k<j$. By hypothesis
the alcove $A'_{i}$ lies in $H^{w+}$. By induction, we know that $A'_k$ also
lies in $H^{w+}$ for any $k=i,\ldots,j-1$ since the only hyperplane 
separating $A'_{k-1}$ and $A'_{k}$ is $H_k\not=H$. Thus $A_{j-1}\in H_j^{w+}$,
which contradicts the hypothesis $A_{j-1}\in H_j^{w-}$.

In an arbitrary minimal gallery $A'_0,\ldots,A'_l$, for every $i$, the alcoves $A'_0$ and 
$A'_{i-1}$ lie in the same connected component of $V-H_i$, 
and the alcoves $A'_{i}$ and $A'_l$ lie in the other connected component. 
If the gallery $A'_0,\ldots,A'_l$ is supposed to be in the $w$-direction,
$A'_{i-1}\in H_i^{w-}$ and $A'_{i}\in H_i^{w+}$ and consequently 
we have $A'_0\in  H_i^{w-}$ and $A'_l\in H_i^{w+}$.
This means for any hyperplane $H$ separating 
$A'_0$ and $A'_l$, $A'_0\in H^{w-}$ and $A'_l\in H^{w+}$.
The same argument now proves that any other minimal gallery going 
from $A'_0$ to $A'_l$ is also in the $w$-direction.
\end{proof}

\begin{defn}
Let $C(A,w)$ be the set of alcoves  $A' \in {\mathcal A}$ such that there exists
a gallery $A'_0, A'_1, \ldots, A'_l$  going from $A$ to $A'$ in the $w$-direction.
\end{defn}  
We call $C(A,w)$ the {\em acute cone from $A$ in the $w$-direction}.

\begin{prop}
We have the inclusion $C_a(A,w)\subset C(A,w)$ for any point $a\in A_0$.
\end{prop}

\begin{proof}
If $A'(a)\in A(a)+w(\bar C_0)$ then for every root $\alpha\in w(R^+)$,
the inequality 
$$\langle \alpha, A(a) \rangle \leq \langle \alpha, A'(a) \rangle$$
holds. Hence, for any hyperplane $H$ separating $A$ and $A'$, 
$A$ lies in the $w$-negative half-space $H^{w-}$ and $A'$ lies
in the $w$-positive half-space $H^{w+}$. As shown in the proof 
of the above lemma, any minimal gallery going from $A$ to $A'$ 
is in the $w$-direction.
\end{proof}

\begin{cor}
For any two alcoves $A,A'$, there exists $w\in W_0$ and a gallery going
from $A$ to $A'$ in the $w$-direction. In other words
$$\calA=\bigcup_{w\in W_0} C(A,w).$$
\end{cor}

\begin{proof}
The union $\bigcup_{w\in W_0} C_a(A,w)$ obviously covers the whole set $\calA$
for any $a\in A_0$.
\end{proof}

\begin{cor}
Let $\mu$ be a dominant coweight. Then $t_{w\mu}(A_0)\in C(A_0,w)$.
\end{cor}

\begin{proof}
If $\mu$ is dominant and lies in $\Check Q$, $t_{w\mu}(A_0)$ belongs to 
$C_a(A_0,w)$ for all $a\in A_0$. If now $\mu$ is dominant but not necessarily 
in $\Check Q$, let us choose a point $a$ in $A_0$ fixed by $\Omega$, for instance
we can take the barycenter of $A_0$. Then
$t_{w\mu}(A_0)$ belongs to $C_a(A_0,w)$ a fortiori to $C(A,w)$.
\end{proof}

There exists another characterization of the acute cone $C(A,w)$ which will be useful later. 

\begin{lemma}
An alcove $A'$ lies in $C(A,w)$ if and only if it is contained in every $w$-positive halfspace
$H^{w+}$ containing $A$:
$$
A'\in \bigcap_{A\in H^{w+}} H^{w+}.
$$
\end{lemma}

\begin{proof}
This is only a matter of rephrasing the proof of Lemma 5.3. In fact, we have seen there
that $A'\in C(A,w)$ if and only if for any wall $H$ separating $A$ and $A'$, we have
$A\in H^{w-}$ and $A'\in H^{w+}$. In other words, $A'\in \bigcap_{A\in H^{w+}} H^{w+}.$
\end{proof}

\section{Proof of Proposition 4.3} 

Suppose $\mu \in X_*$ is dominant and $\mu \in W_{\rm aff}\tau$ for $\tau \in \Omega$, 
the isotropy group of the base alcove.  Given an element $x\in W_{\rm aff}\tau$ which 
is $\mu$-strongly permissible, we want to prove there  
exists $w\in W_0$ such that $x\leq t_{w\mu}$. According to the Corollary 5.6,
there exists $w\in W_0$ such that $x(A_0) \in C(A_0,w)$. According to Corollary 5.7, 
$t_{w\mu}(A_0)$ lies in the same acute cone $C(A_0,w)$. We will prove $x \leq t_{w\mu}$
for this $w$.

Let us write $x=x_1\tau$ and $t_{w\mu} =x_2\tau$ with $x_1, x_2\in W_{\rm aff}$
and $\tau\in\Omega$. To prove Proposition 2 we thus have 
to prove the following statement.

\begin{prop}
Let $x_1, x_2$ be two elements of $W_{\rm aff}$ and $w\in W_0$ such that
\begin{itemize}
\item for every vertex $v$ of $A_0$, $x_1(v)$ lies in $B(x_2(v),w)$;
\item the two alcoves $x_1(A_0)$ and $x_2(A_0)$ lie in $C(A_0,w)$.
\end{itemize}
Then the inequality $x_1\leq x_2$ holds. 
\end{prop}

First, we prove that the statement does not change when we move the base alcove
in the $w$-opposite direction. By this operation we can then move the alcoves $x(A_0)$ and 
$t_{w\mu} (A_0)$ far from the walls. As we will see, the proof is much easier
in that case, compared to the proof of Kottwitz-Rapoport in the minuscule case,
which deals with alcoves very close to the origin. 

\begin{lemma}
Let $x_1,x_2\in W_{\rm aff}$ be such that the alcoves $x_1(A_0)$ and $x_2(A_0)$ lie
in $C(A_0,w)$. Let $y\in W_{\rm aff}$ be such that $A'_0\in C(A_0,ww_0)$ where 
$y(A'_0)=A_0$.  Then $x_1\leq x_2$ if and only if $yx_1\leq yx_2$.
\end{lemma}

Note that the condition $A'_0\in C(A_0,ww_0)$ means there exists a gallery 
going from $A'_0$ to $A_0$ in the $w$-direction.

\begin{proof}
According to Bernstein-Gelfand-Gelfand \cite{BGG} (see also \cite{Deod1})
every Coxeter group satisfies the following property :
let $x_1,x_2\in W_{\rm aff}$ and $s\in S_{\rm aff}$ be a simple reflection, then
$$\begin{array}{rcl}
x_1\leq x_2 & \Leftrightarrow& sx_1\leq sx_2
\end{array}$$ 
whenever $l(sx_1)>l(x_1)$ and $l(sx_2)>l(x_2)$.  
It follows easily by induction on $l(y)$ that $x_1 \leq x_2$ is equivalent to $yx_1 \leq yx_2$
whenever $l(yx_i) = l(y) + l(x_i)$ for $i=1,2$.

Therefore what we really have to prove is that the lengths add
$$\begin{array}{rcl}
l(yx_1)=l(y)+l(x_1)& \ \ & l(yx_2)=l(y)+l(x_2).
\end{array}$$
Of course, we only need to prove the first of these equalities.

Let us choose a minimal gallery $A'_0,\ldots,A'_l$ 
going from $A'_0$ to $A_0$, with $l=l(y)$, and $A_0,\ldots,A_r$
a minimal gallery going from $A_0$ to $x_1(A_0)$, with $r=l(x_1)$. 
According to Lemma 5.3,
any minimal galleries going from $A'_0$ to $A_0$ or going from $A_0$
to $x_1 (A_0)$ are automatically in the $w$-direction. Obviously, the concatenation of two 
galleries in the $w$-direction is still in the $w$-direction. By Lemma 5.3 again,
the concatenated gallery $A'_0, \dots, A'_l = A_0, \dots, A_r$ is minimal.

Let us prove that the concatenated gallery 
$A'_0,\ldots,A_0,\ldots,A_r=x_1yA'_0$ is minimal only if 
the equality  $l(yx_1)=l+r$ holds.

Let $S'_{\rm aff}=\lbrace y^{-1}sy\mid s\in S_{\rm aff}\rbrace$ be
the set of reflections by walls of $A'_0$; obviously $(W_{\rm aff},S'_{\rm aff})$
is also a Coxeter system. Now, by viewing $A'_0$ as base alcove,
the gallery $A'_0,\ldots,A_0,\ldots,A_r=x_1yA'_0$ gives rise to a reduced 
expression 
$$x_1y=s'_1\ldots  s'_{l+r}$$
with $s'_1,\ldots,s'_{l+r}$  in $S'_{\rm aff}$. Thus the length 
of $x_1y$ in the Coxeter system $(W_{\rm aff},S'_{\rm aff})$ is 
equal to $l+r$.
Consequently, the length of 
$$yx_1=(ys'_1 y^{-1})\ldots (ys'_{l+r} y^{-1}) $$
in the Coxeter system $(W_{\rm aff},S_{\rm aff})$ is equal to $l+r$.
\end{proof}

Another important ingredient in the proof is a lemma due to Deodhar \cite{Deod1},
which is available for any Coxeter system, in particular for the Coxeter system 
$(W_{\rm aff},S_{\rm aff})$. Let $J$ be a subset of $S_{\rm aff}$ ; 
let $W_J$ be the subgroup  of $W_{\rm aff}$ generated by $J$. It is known 
that there exists a subset $W^J$ of $W_{\rm aff}$  such that any element 
$x\in W_{\rm aff}$ admits a {\em unique decomposition}
$x=x^J x_J$ with $x^J\in W^J$ and $x_J\in W_J$ such that $l(x)=l(x^J)+l(x_J)$.

\begin{lemma} [{\bf Deodhar [3]}]
Let $\mathcal J$ be a set whose elements are subsets of $S_{\rm aff}$ such that
$\bigcap_{J\in\mathcal J}J=\emptyset$. Let $x,y\in W_{\rm aff}$ and 
$x=x^Jx_J$ and $y=y^Jy_J$ be the unique decompositions of $x$ and 
$y$ for $J\in\mathcal J$. Then $x\leq y$ if and only if $x^J\leq y^J$ for every 
$J\in\mathcal J$. 
\end{lemma}

Elements of $S_{\rm aff}$ correspond to walls of the base alcove $A_0$. 
For any vertex $a$ of $A_0$, let $J_a$ denote the subset of elements of $S_{\rm aff}$
corresponding to walls containing $a$. Obviously, the intersection of such 
subsets $J_a$ for all vertices $a$ of the base alcove is empty. Thus, the lemma of Deodhar 
applies.

Let us fix a vertex $a$ of the base alcove $A_0$ and let $J$ denote $J_a$.
It is useful to visualize the decomposition $x=x^{J}x_{J}$ in terms
of alcoves. Let $x$ be an element of
$W_{\rm aff}$ and $A=x(A_0)$ the corresponding alcove. The coset $xW_J$ 
 is in bijection with the set of alcoves sharing with $A$ the
vertex $x(a)$. The element $x^J$ corresponds to the minimal element
among those alcoves.

\begin{lemma}
Let $a$ be a vertex of the base alcove $A_0$ and let $J$ denote $J_a$.
 Let $x_1,x_2$ be elements of 
$W_{\rm aff}$ such that $x_1(a)\in B(x_2(a),w)$. There exists an alcove $A$ 
such that for every $y\in W_{\rm aff}$ satisfying $y^{-1}A_0\in C(A,ww_0)$
we have $(yx_1)^{J}\leq (yx_2)^{J}$. 
\end{lemma}

In more intuitive terms, the inequality $(yx_1)^{J}\leq (yx_2)^{J}$
holds for the alcove $y^{-1}A_0$ far enough in $w$-opposite direction.

\begin{proof}
Let us denote $v_1=x_1(a)$ and $v_2=x_2(a)$.
We can easily reduce to the case $v_1=s v_2$, where $s$ is some 
reflection $s_{\alpha,k}$ and $v_1\in v_2+w(B_0)$. 
In other words $v_1$ and $v_2$ are symmetric 
with respect to the hyperplane $H=H_{\alpha,k}$ ; $v_1$ lies in the $w$-negative 
half-space $H^{w-}$ ; $v_2$ lies in the $w$-positive half-space
$H^{w+}$. We can also suppose $v_1\not=v_2$, thus $v_1$ and $v_2$ lie in the interiors 
of the corresponding half-spaces determined by $H$. 

Let $A$ be an alcove in the $w$-negative half-space
$H^{w-}$. The whole acute cone $C(A,ww_0)$ belongs then to 
$H^{w-}$.

Let $A'_0=y^{-1}(A_0)$ be an alcove lying in $C(A,ww_0)$, a fortiori in the 
$w$-negative half-space $H^{w-}$. Just as in the proof of the last lemma,
it will be convenient to consider $A'_0$ as base alcove. Let $S'_{\rm aff}$
denote the set of reflections by walls of $A'_0$. Let $a'=y^{-1}(a)$ and let $J'$
denote the subset of reflections by the walls of $A'_0$ containing the vertex $a'$.
Since the conjugation by $y^{-1}$ makes the Coxeter $(W_{\rm aff}, S_{\rm aff})$ 
isomorphic to the Coxeter system $(W_{\rm aff}, S'_{\rm aff})$, the inequality 
we want to prove in the former system 
$$(yx_1)^J\leq (yx_2)^J$$
is equivalent to the inequality we are going to prove
$$(x_1y)^{J'}\leq (x_2y)^{J'}$$
in the latter system. In the remainder of the proof, all inequalities are understood 
to be relative to the Coxeter system $(W_{\rm aff}, S'_{\rm aff})$.

Obviously, $(x_1y)^{J'}(a')=x_1y(a')=v_1$ and  $(x_2y)^{J'}(a')=x_2y(a')=v_2$.
The alcove $(x_1y)^{J'}(A'_0)$, resp. $(x_2y)^{J'}(A'_0)$, is minimal among the
alcoves sharing the vertex $v_1$, resp. $v_2$. 
Like $v_1$, the alcove $(x_1y)^{J'}(A'_0)$ lies in the $w$-negative half-space 
$H^{w-}$.  Like $v_2$, the alcove $(x_2y)^{J'}(A'_0)$ lies in the $w$-positive 
half-space $H^{w+}$. 

Let $A'_0,A'_1,\ldots,A'_l$ 
be a minimal gallery going from $A'_0$ to $A'_l=(x_2y)^{J'} A'_0$.
Since $A'_0\in H^{w-}$ and $A'_l\in H^{w+}$, there
exists an integer $j$ such that $A'_{j-1}$ and $A'_{j}$ share the wall $H$.
The minimal gallery $A'_0,A'_1,\ldots,A'_l$ corresponds to a reduced expression
$$(x_2y)^{J'} =s'_1 s'_2 \ldots s'_l$$
with $s'_1,\ldots,s'_l\in S'_{\rm aff}$.
By removing the reflection $s'_j$ from the reduced expression, we get 
$$s'_1\ldots \hat s'_j\ldots s'_l=s (x_2y)^{J'}.$$
Therefore $s(x_2y)^{J'}\leq (x_2y)^{J'}$. 

But the alcove $s(x_2y)^{J'}A'_0$ contains the vertex $s(x_2y)^{J'}(a')=v_1$.
Thus by minimality, we know $(x_1y)^{J'}\leq s (x_2y)^{J'}$,
therefore $(x_1y)^{J'}\leq (x_2y)^{J'}$.
\end{proof}

\noindent{\em Proof of Proposition 6.1.}
The proposition follows from Lemmas 6.2, 6.3 and 6.4 as follows.
By hypothesis we know that $x_1(A_0), \,\, x_2(A_0) \in C(A_0,w)$ and that, for every vertex $a$ of the base alcove, $x_1(a)\in B(x_2(a),w)$.  By 6.4, for every vertex $a \in A_0$ there exists an alcove $A_a$ such that, if $y \in W_{\rm aff}$ satisfies $y^{-1}(A_0) \in C(A_a, ww_0)$, then $(yx_1)^{J_a} \leq (yx_2)^{J_a}$.  Now choose $y$ such that
$$
y^{-1}(A_0) \in C(A_0,ww_0) \cap \bigcap_{a} C(A_a,ww_0).
$$
(Note that the intersection of finitely many acute cones in the same direction is nonempty.)  Then $(yx_1)^{J_a} \leq (yx_2)^{J_a}$ for every $a$ and thus by Lemma 6.3, $yx_1 \leq yx_2$.  But then Lemma 6.2 implies $x_1 \leq x_2$.

\hfill $\square$

\section{Deodhar's counter-examples}

In \cite{Deod2}
Deodhar considers the following question: which finite Weyl groups 
$W_0$ have the property that  $w \leq w'$ in the Bruhat order if and only if 
$w(\lambda) - w'(\lambda)$ is a sum of positive 
coroots for every dominant coweight $\lambda$?  
Surprisingly enough, at least for us, he proves that $W_0$ 
has this property if and only if the irreducible components of the 
associated root system $R$ are of type $A_{n-1}$ or of rank $\,\, \leq 3$.  
For part of his proof, 
for an irreducible root system of rank $\geq 4$ and not of type $A_{n-1}$, 
he proves the following statement by giving explicit examples.

\begin{prop}[{\bf Deodhar [4]}] 
Let $R$ be an irreducible root system of rank $\geq 4$ and not of type $A_{n-1}$.  Then there exist elements $w,w'\in W_0, \,\, w \neq w'$ such that
\begin{itemize}
\item for every dominant coweight $\lambda$, $w(\lambda) - w'(\lambda)$ 
is a sum of positive coroots;
\item $l(w)=l(w')$.
\end{itemize}

\end{prop}

We derive the following statement about $\mu$-admissible and 
$\mu$-permissible sets, thus proving Theorem 3.

\begin{prop}
For $R$ an irreducible root system of rank $\geq 4$ and not of type $A_{n-1}$,
we have $\Adm(\mu)\not=
\Perm(\mu)$ for every sufficiently regular dominant coweight $\mu$.
\end{prop}

\begin{proof}
The idea of the proof is the remark, due to Kottwitz and Rapoport, that nearby the extreme 
elements $t_{w\mu}$, the picture looks like the Bruhat order in the finite Weyl 
group $W_0$.  (This is proved in a precise form in Lemma 7.5 below.)

Let $\mu$ be a {\em regular} dominant coweight. For simplicity, suppose
$\mu\in \Check Q$. Let $w,w'\in W_0$ be chosen as in Proposition 7.1. Let us consider
the element
$$x=t_{w^{-1}\mu}w^{-1}w'.$$
We prove $x$ is $\mu$-permissible but not $\mu$-admissible, if $\mu$ is sufficiently regular.

Since $w \neq w'$, we know $x\not= t_{\mu'}$ for any $\mu'\in W_0(\mu)$.
Therefore to prove $x$ is not $\mu$-admissible, it is enough to show $l(x)=l(t_\mu)$.
Since $\mu$ is regular dominant, by a theorem of Iwahori-Masumoto \cite{IM} 
the maximal length element in the coset $t_{w^{-1}\mu}W_0$ is $t_{w^{-1}\mu}w^{-1}$.
Moreover we have 
$$\begin{array}{rcl}
l(x)&=&l(t_{w^{-1}\mu}w^{-1})-l(w')\\
&=& l(t_{w^{-1}\mu}w^{-1})-l(w)\\
&=& l(t_{w^{-1}\mu}),
\end{array}$$
so $l(x)$ is equal to $l(t_\mu)$.

Let us prove that $x$ is $\mu$-permissible.  Let $a$ be a vertex of $A_0$. 
By construction of the elements $w$ and $w'$, the difference
$$t_{w^{-1}\mu}(a)-t_{w^{-1}\mu}w^{-1}w' (a)=a-w^{-1}w' (a)=w^{-1}(w(a)-w'(a))$$
lies in $-w^{-1}(B_0)$, since $a$ is dominant.

For $\mu$ sufficiently regular dominant, we have
$$t_{w^{-1}\mu}w^{-1}w' (a)-a\in w^{-1}(\bar C_0).$$
(This relation can be used as the {\em definition} of ``sufficiently regular''.)  According to the following well-known statement, $x$ is $\mu$-permissible.
\end{proof}

\begin{lemma} The equality
$$w\bar C_0\cap (w\mu+ wB_0)=w\bar C_0\cap  {\rm Conv}(\mu)$$
holds for any $\mu$ dominant and $w\in W_0$.
\end{lemma}

\noindent In some sense, this lemma was the starting point of this work.

\begin{Remark}
The following lemma, which is not used in the sequel, makes precise the remark of Kottwitz-Rapoport which inspired Proposition 7.2 above.
\end{Remark}

\begin{lemma}
Let $R$ be any irreducible root system.  Suppose $\lambda \in X_*$, and let $t_\lambda w_\lambda$ be the unique element of minimal length in the coset $t_\lambda W_0$.  Suppose $w_1, \, w_2 \in W_0$.  Then $t_\lambda w_\lambda w_1 \leq t_\lambda w_\lambda w_2$ if and only $w_1 \leq w_2$.  
\end{lemma}

\begin{proof}
By a formula of Iwahori-Matsumoto (loc. cit.), $l(t_\lambda w_\lambda w_i) = l(t_\lambda w_\lambda ) + l (w_i)$ for $i=1,2$.  The lemma is a consequence of these equalities, as explained already in the proof of Lemma 6.2.
\end{proof}

The geometric meaning of the lemma can be stated as follows.

\begin{cor}
Suppose $x, \, y \in W_{\rm aff}\tau$ are such that the alcoves $x(A_0)$ and $y(A_0)$ share a vertex $v = t_\lambda(0)$, where $t_\lambda \in W_{\rm aff}\tau$.   Then $x < y$ if and only if there exists a sequence of reflections $s_{H_1}, \dots, s_{H_n}$ such that
$$
x < s_{H_1}x < s_{H_2}s_{H_1}x < \cdots < s_{H_n}\dots s_{H_1}x = y,
$$
where every hyperplane $H_i,  \,\ i=1,\dots,n$, contains the vertex $v$.
\end{cor}

\section{Proof of Deodhar's theorem for $A_{n-1}$}

The goal of this section is to present a short proof of the following theorem of Deodhar \cite{Deod2} characterizing the Bruhat order on $W_0 = S_{n-1}$.  As in the proof of Theorem 1 we rely on Proposition 6.1 as the key ingredient.    

\begin{prop} [{\bf Deodhar [4]}]   Let $R$ be a root system of type $A_{n-1}$.  Let $w, w' \in W_0$.  Then $w' \leq w$ in the Bruhat order on $W_0$ if and only if for every dominant coweight $\lambda$, $w'(\lambda) - w(\lambda)$ is a sum of positive coroots.
\end{prop}

\begin{proof}
We need to show that if $w'(\lambda) -  w(\lambda)$ is a sum of positive coroots for all dominant coweights $\lambda$, then $w' \leq w$, the reverse implication being a general fact for all root systems which is easily proved.  Recall that for root systems of type $A_{n-1}$ every vertex $v$ of the base alcove is either a dominant fundamental coweight or zero.  Let $w_0$ denote the longest element of $W_0$.  
  
The assumption on $w$ and $w'$ implies that, for every vertex $v \in A_0$,
$$w'(v) \in W_{\rm aff}(w(v)) \cap (w(v) + w_0(B_0))$$
and thus
$$w'(v) \in B(w(v),w_0),$$
by Lemma 4.4.
By Lemma 8.2 below, $w(A_0)$ and $w'(A_0)$ both belong to $C(A_0,w_0)$.  Therefore by Proposition 6.1, $w' \leq w$, as desired. 

\end{proof}

We have used a special case ($\mu = 0$) of the following general lemma.  Intuitively, it says that all the alcoves ``on the boundary'' of the antidominant Weyl chamber are contained in the antidominant acute cone of alcoves.  

\begin{lemma}
Let $R$ be an irreducible root system with Weyl group $W_0$.  Let $w_0$ denote the longest element of $W_0$.  Then for any dominant coweight $\mu$ and any $w \in W_0$ we have
$$
t_{w_0 \mu} w(A_0) \in C(A_0,w_0).
$$
\end{lemma}

\begin{proof}
Since the concatenation of two galleries in the $w_0$-direction is still in the $w_0$-direction, it suffices to prove that $t_{w_0 \mu}(A_0) \in C(A_0,w_0)$ and $t_{w_0 \mu} w(A_0) \in C(t_{w_0 \mu}(A_0),w_0)$.  The first statement results from Corollary 5.7, and since the translation of a gallery in the $w_0$-direction is still in the $w_0$-direction, the second is equivalent to $w(A_0) \in C(A_0,w_0)$.

Choose a reduced expression $w = s_1s_2\cdots s_r$, where for each $i = 1, \dots r$, $s_i = s_{\alpha_i}$  is the reflection corresponding to a simple root $\alpha_i$.  We claim that the gallery $A_0,s_1(A_0), \dots, s_1\cdots s_r(A_0)$ is in the $w_0$-direction.  
Fix $i$ and let $H_i = H_{s_1\cdots s_{i-1}(\alpha_i)}$ denote the 
hyperplane separating $s_1\cdots s_{i-1}(A_0)$ and $s_1\cdots s_i(A_0)$.  We need to show that
$$
s_1\cdots s_{i-1}(A_0) \subset H^{w_0-}_i 
$$
and 
$$
s_1\cdots s_i(A_0) \subset H^{w_0+}_i.
$$
But since $s_1\cdots s_{i-1}(\alpha_i) > 0$, we have
$$
H^{w_0-}_i = \lbrace x \in V ~| ~ \langle s_1\cdots s_{i-1}(\alpha_i),x \rangle > 0 \rbrace
$$ 
and 
$$
H^{w_0+}_i = \lbrace x \in V ~ | ~ \langle s_1\cdots s_{i-1}(\alpha_i),x \rangle < 0 \rbrace,
$$
yielding the desired inclusions.

\end{proof}

\begin{Remark}
If $R$ is any root system such that
$$
B(v,w_0) = W_{\rm aff}(v) \cap (v + w_0(B_0))
$$
for every vertex $v$, then as in the proof of Proposition 8.1. we see that $W_0$ enjoys the property that $w' \leq w$ in the Bruhat order if and only if $w'(\lambda) - w(\lambda)$ is a sum of positive coroots for every dominant coweight $\lambda$.  (Indeed, letting 
$\lambda$ range over fundamental coweights and recalling that every vertex of $A_0$ is a 
positive scalar multiple of a fundamental coweight, we see that the condition on the differences $w'(\lambda) - w(\lambda)$ ensures that 
$w'(v) \in W_{\rm aff}(w(v)) \cap (w(v) + w_0(B_0))$ for every vertex $v \in A_0$, and 
the rest of the proof goes over word-for-word.)  Consequently, Deodhar's theorem (see Proposition 7.1) implies that $R$ is of type $A_{n-1}$ or of rank $\leq 3$ (cf. Remark 4.5). 
\end{Remark}

\section{Admissible sets and automorphisms of root systems}

In this section we study the behavior of admissible sets relative to an 
automorphism $\Theta$ of the underlying root system.  This will yield information about the 
admissible and permissible sets in the extended affine Weyl group for ${\rm GSp}_{2n}$, to 
be discussed at the end of this section and in section 10.  We begin by recalling some facts about automorphisms of root 
systems; the main reference is Steinberg's article \cite{Steinberg}, especially $\S 1.30-1.33$. 

Let $R = (X^*,X_*,R,\Check{R},\Pi)$ be an irreducible and reduced (based) root system, as in 
section 2.  We may equip $V = X_* \otimes \RR$ with a $W_0$-invariant inner product $(~,~)$.  
By assumption $X^* = {\rm Hom}(X_*,\ZZ)$, and we let $\langle ~, ~ \rangle : X^* \times X_* 
\rightarrow \ZZ$ denote the canonical pairing.  A root $\beta \in R$ is a linear functional on 
$V$.  Using the isomorphism of the vector space $V$ with its dual defined with the help of 
$(~,~)$, we have the identification $\beta = 2\Check{\beta}/(\Check{\beta},\Check{\beta})$; in 
this way $\beta$ may be regarded as an element of $V$.

Let $\Theta$ be an automorphism of $R$, in the sense of \cite{Steinberg}.  This means that 
$\Theta$ is a linear automorphism of $(V,(~,~))$ leaving stable the sets $X^*$, $R$, and $\Pi$, 
when these are regarded as functionals on $V$.  Hence, $\Theta$ leaves stable the subsets 
$X_*$ and $\Check{R}$ of $V$ as well.   Clearly $\Theta$ induces automorphisms of the 
groups $\widetilde{W}$, $W_{\rm aff}$, $W_0$, and $\Check{Q}$; in particular we may consider 
the fixed-point subgroups $W^\Theta_{\rm aff}$, $W^\Theta_0$, and $\Check{Q}^\Theta$.

The group $W^\Theta_0$ is the Weyl group of a root system $R^{[\Theta]}$, which is defined 
as follows.  Define ${\mathcal Z} = \lbrace x \in X_* ~ | ~ \langle \alpha , x \rangle = 0, \,\, 
\mbox{for all} \,\, \alpha \in R \rbrace$.  It is known that ${\mathcal Z} \cap \Check{Q} = \lbrace 0 
\rbrace$ and ${\mathcal Z} + \Check{Q}$ has finite index in $X_*$ (see Lemma 1.2 of \cite{Springer}). 
Let $X^{[\Theta]}_* = \lbrace x \in X_* ~ | ~ x - \Theta(x) \in {\mathcal Z} \rbrace$, and let $V^{[\Theta]} 
= X^{[\Theta]}_* \otimes \RR$.  For any $\beta \in R$, let $\bar{\Theta}(\beta)$ denote the average 
of the $\Theta$-orbit of $\beta$.  Then according to \cite{Steinberg} $\S 1.33$, $R^{[\Theta]}$ is 
the subset of $V^{[\Theta]}$ consisting of all elements $\bar{\Theta}(\beta)$, $\beta \in R$, except
 those which are smaller multiples of others.  Moreover, for any $\Theta$-orbit $\pi \subset \Pi$, 
let $V_\pi$ denote the real vector space generated by $\pi$ and let $\tilde{\alpha}_\pi$ be any 
highest root of the (possibly reducible) root system $R \cap V_\pi$.  Then the set $\Pi^{[\Theta]}$ 
of simple roots in $R^{[\Theta]}$ consists of the elements 
$\alpha_\pi := \bar{\Theta}(\tilde{\alpha}_\pi)$ as $\pi$ ranges over all $\Theta$-orbits in $\Pi$. 
(The highest roots of $R \cap V_\pi$ form a single $\Theta$-orbit, so the choice of $\tilde{\alpha}_\pi$ 
for each orbit $\pi$ is immaterial.)  By loc. cit. $\S 1.32-1.33$, $(V^{[\Theta]}, R^{[\Theta]}$, 
$\Pi^{[\Theta]})$ is a root system in the sense of loc. cit., $\S 1.1-1.6$.         

For any $\Theta$-orbit $\pi \subset \Pi$, let $s_\pi$ denote the longest element of the Weyl 
group $W_\pi$ generated by $\lbrace s_\alpha ~ | ~ \alpha \in \pi \rbrace$.  The restriction of 
$s_\pi$ to $V^{[\Theta]}$ is the reflection $s_{\alpha_\pi}$ corresponding to $\alpha_\pi$.  
In the sequel we will abuse notation and write $s_{\alpha_\pi}$ instead of $s_\pi$.

Let $X^{*[\Theta]} = {\rm Hom}(X^{[\Theta]}_*,\ZZ)$.  Regarding $R^{[\Theta]}$ as a set 
of linear functionals on $V^{[\Theta]}$, it follows easily that $R^{[\Theta]} \subset 
X^{*[\Theta]}$.  Regarding $R^{[\Theta]}$ as a subset of $V^{[\Theta]}$, we define another subset 
$\Check{R}^{[\Theta]} = \lbrace 2\beta/(\beta,\beta) ~ | ~ \beta \in R^{[\Theta]} \rbrace$.  
The following statement is certainly well-known, but lacking a convenient reference, we 
provide some details.

\begin{prop}
$R^{[\Theta]} = (X^{*[\Theta]},X^{[\Theta]}_*, R^{[\Theta]},\Check{R}^{[\Theta]}, \Pi^{[\Theta]})$ 
is a reduced and irreducible (based) root system in the vector space $V^{[\Theta]}$.  Its affine 
Weyl group is $W_{\rm aff}^\Theta = \Check{Q}^\Theta \rtimes W^\Theta_0$.
\end{prop}

\begin{proof}
In light of loc. cit. $\S 1.32$, to prove $R^{[\Theta]}$ is a based root system it remains only 
to show that, (i) every element of $R^{[\Theta]}$ is an {\em integral} linear combination of 
elements of $\Pi^{[\Theta]}$, and (ii) $\Check{R}^{[\Theta]} \subset X^{[\Theta]}_*$.  For 
statement (i), consider the following property of a $\Theta$-orbit $\pi \subset \Pi$:
$$
\noindent (*) \phantom{....}  \mbox{{\em The elements of $\pi$ are pairwise orthogonal}.}
$$
If $(*)$ holds for every $\pi$, then each $\alpha \in \pi$ is a highest root and thus 
$\bar{\Theta}(\alpha) = \alpha_\pi$.  Therefore (i) follows from the analogous property 
of the root system $R = (X^*, X_*, R, \Check{R}, \Pi)$.  By examining automorphisms of 
simple Dynkin diagrams, we see that $(*)$ holds for every 
$\pi$ unless $R$ is of type $A_{2n}$ (which we assume realized in 
$V = \RR^{2n+1}$), $\Theta(x_1,\dots,x_{2n+1}) = (-x_{2n+1},\dots,-x_1)$, and 
$\pi = \lbrace e_n - e_{n+1}, e_{n+1} - e_{n+2} \rbrace$; therefore, $\alpha_\pi = e_n - e_{n+2}$. 
In this case one can verify by direct calculation that (i) holds 
(and in fact $R^{[\Theta]}$ is of type $C_n$). 
We note that (i) is equivalent to the statement that $2(\alpha,\beta)/(\beta,\beta) \in \ZZ$ for 
every $\alpha, \beta \in R^{[\Theta]}$. 

Next, one can show, following loc. cit. $\S 1.32$, that $W^{\Theta}_{\rm aff}$ is the group of 
affine transformations of $V^{[\Theta]}$ generated by the reflections through the hyperplanes 
$\beta + k = 0$, for $\beta \in R^{[\Theta]}$ and $k \in \ZZ$.  It follows that $W^{\Theta}_{\rm aff}$ 
is the affine Weyl group for the root system $R^{[\Theta]}$ in $V^{[\Theta]}$, and therefore 
the fixed-point group $\Check{Q}^\Theta$ is precisely the corresponding coroot lattice.  But 
then $\Check{R}^{[\Theta]} \subset \Check{Q}^\Theta \subset X^{[\Theta]}_*$, and (ii) is proved.

Finally, if $\tilde{\alpha}$ is the highest root of $R$, then $\Theta(\tilde{\alpha})$ is also highest 
and thus $\Theta(\tilde{\alpha}) = \tilde{\alpha}$ is the unique highest root of $R^{[\Theta]}$.  
Therefore $R^{[\Theta]}$ is irreducible. 
\end{proof}

Viewing roots as linear functionals on $V$, we have the inclusion
$$
R^{[\Theta]} \subset R|_{V^{[\Theta]}},
$$
since $\alpha|_{V^{[\Theta]}}$ can be identified with $\bar{\Theta}(\alpha)$, 
for any $\alpha \in R$. Although the opposite inclusion does not always hold, 
every element of $R|_{V^{[\Theta]}}$ is of the form $c \beta$ for some 
$\beta \in R^{[\Theta]}$ and $c \in [0,1]$. In fact $c\in\{{1\over 2},1\}$, as is 
shown by the 
following lemma which  is implicit in Steinberg's book \cite{Steinberg}; 
we provide a proof since we could not locate a reference.

\begin{lemma}
Given a root $\alpha \in R$, $\bar{\Theta}(\alpha)$ is either a root $\alpha'$ or a half-root 
$\frac{1}{2}\alpha'$, for $\alpha' \in R^{[\Theta]}$.  In other words
$$
\bar{\Theta}(R) \subset 
W^\Theta_0 \Pi^{[\Theta]} \cup \frac{1}{2}W^\Theta_0 \Pi^{[\Theta]}.
$$
Conversely, if $\alpha'$ is any root in $R^{[\Theta]}$, 
there is at least one root $\alpha \in R$ such that $\bar{\Theta}(\alpha) = \alpha'$.  
Moreover $\alpha'$ is a positive root in $R^{[\Theta]}$ if and only if $\alpha$ is a 
positive root in $R$.
\end{lemma}
 
\begin{proof}
First note that if $\alpha \in V_\pi \cap R$, then $\bar{\Theta}(\alpha)$ is either $\alpha_\pi$ 
or $\frac{1}{2} \alpha_\pi$ (the latter occurring only in the case $R$ is of type $A_{2n}$ and 
$\pi = \lbrace e_{n} - e_{n+1} , e_{n+1} - e_{n+2} \rbrace$; cf. the proof of Prop. 9.1).  Thus it 
is enough to prove the following statement.

\medskip
$(\dagger)$ \hspace{.2in} Given $\beta \in R^+$, there exist $\pi \subset \Pi$ and 
$w \in W^\Theta_0$ such that $w\beta \in V_\pi \cap R$.    
\medskip 

\noindent For any $\beta \in R^+$, we have $\bar{\Theta}(\beta) = \sum_\pi c_\pi \alpha_\pi$, 
where $c_\pi \geq 0$ for every $\Theta$-orbit $\pi \subset \Pi$.  We will prove $(\dagger)$ by 
induction on the quantity $\Sigma \bar{\Theta}(\beta) = \sum_\pi c_\pi$.

Fix a positive root $\beta$; we may assume that $\beta$ is not supported in any $\pi \subset \Pi$.  
Since the pairing $(~,~)$ is positive definite, we have $(\bar{\Theta}(\beta),\bar{\Theta}(\beta)) > 0$ 
and thus there is a $\Theta$-orbit $\pi$ such that $(\bar{\Theta}(\beta),\alpha_\pi) > 0$. 
By \cite{Steinberg} $\S 1.15$, $s_{\alpha_\pi}$ permutes the positive roots not supported on $\pi$, 
and thus $s_{\alpha_\pi}\beta$ is positive.  On the other hand

$$
\bar{\Theta}(s_{\alpha_\pi}\beta) = s_{\alpha_\pi}\bar{\Theta}(\beta) = 
\bar{\Theta}(\beta) - \frac{2(\bar{\Theta}(\beta),\alpha_\pi)}{(\alpha_\pi,\alpha_\pi)}\alpha_\pi,
$$
and thus $\Sigma \bar{\Theta}(s_{\alpha_\pi}\beta) <  \Sigma \bar{\Theta}(\beta)$.  
By induction $s_{\alpha_\pi}\bar{\Theta}(\beta) \in W^\Theta_0\Pi^{[\Theta]} \cup 
\frac{1}{2}W^\Theta_0\Pi^{[\Theta]}$, whence the first statement follows.
The converse statement, as well as the one concerning positivity, is obvious.  
\end{proof}
 
Given a fixed affine root $(\alpha,k)$ for $R$, let $k' = k$ if $\bar{\Theta}(\alpha) = \alpha'$ and 
let $k' = 2k$ if $\bar{\Theta}(\alpha) = \frac{1}{2}\alpha'$.  
Given a hyperplane $H_{\alpha,k}$ for $R$ in $V$ we have
$$
H_{\alpha,k} \cap V^{[\Theta]} = H_{\alpha',k'},
$$
where $\alpha'$ and $k'$ are defined as above.  Moreover, every hyperplane for 
$R^{[\Theta]}$ in $V^{[\Theta]}$ is of this form.  This has the following consequence
for the set of alcoves in $V^{[\Theta]}$ that arise from the root system $R^{[\Theta]}$.

\begin{prop}
Given an alcove $A$ in $V$, the set $A \cap V^{[\Theta]}$ is either empty or is an alcove 
$A'$ in $V^{[\Theta]}$; moreover every alcove $A'$ in $V^{[\Theta]}$ arises as 
$A' = A \cap V^{[\Theta]}$ for a uniquely determined alcove $A$ in $V$. 

If $A_0$ (resp. $A'_0$) is the base alcove in $V$ (resp. in $V^{[\Theta]}$), we have 
$A'_0=A_0\cap V^{[\Theta]}$. 
\end{prop}

\begin{proof} If $x,y\in V^{[\Theta]}$ are in the same alcove $A'$, they belong to the same
connected component of $V-\bigcup H_{\alpha,k}$; therefore $x,y\in A$ for some alcove
$A$ in $V$. If $x,y\in V^{[\Theta]}$ belong to two different alcoves in $V^{[\Theta]}$, there exists 
$(\alpha',k')\in R^{[\Theta]}\times\ZZ$ such that $\langle \alpha',x \rangle <k'$ and $\langle \alpha',y \rangle >k'$. Let $(\alpha,k)$
be a corresponding affine root in $R$ (i.e., $\bar{\Theta}(\alpha) = \alpha'$ and $k = k'$). We have $\langle \alpha,x \rangle <k$ and $\langle \alpha,y \rangle > k$, thus $x$ and
$y$ belong to different alcoves in $V$.

An element $x\in A'_0$ satisfies the inequalities $\langle \bar\Theta(\alpha),x \rangle >0$ for all positive roots
$\alpha$ of $R$, and the inequality $\langle \bar\Theta(\tilde\alpha),x \rangle <1$ for the highest root 
$\tilde\alpha$ of $R$. Since $x\in V^{[\Theta]}$ these inequalities are equivalent with 
$\langle \alpha,x \rangle >0$ for all positive roots $\alpha$ of $R$, and  $\langle \tilde\alpha,x \rangle <1$. Therefore 
$x\in A_0$.
\end{proof}

\begin{lemma}
For any hyperplane $H_{\alpha,k}$ for $R$ in $V$, the equality
$$
H^{w+}_{\alpha,k} \cap V^{[\Theta]} = H^{w+}_{\alpha',k'}
$$
holds for any $w \in W^\Theta_0$.
\end{lemma}

\begin{proof}
This follows easily from the definitions, noting that $\alpha \in w(R^+)$ if and  only if 
$\alpha' \in w(R^{[\Theta]+})$ under the assumption $w\in W_0^\Theta$.
\end{proof}

It is now easy to prove an inheritance property for acute cones of alcoves.
Let $C(A_0,w) \cap V^{[\Theta]}$ denote the set of alcoves in $V^{[\Theta]}$ of the form 
$A \cap V^{[\Theta]}$ such that $A \in C(A_0,w)$.

\begin{prop}
Let $w \in W^\Theta_0$.  Then
$$
C(A^{[\Theta]}_0,w) = C(A_0,w) \cap V^{[\Theta]}.
$$
\end{prop}

\begin{proof}
This follows easily from Lemma 9.4 and the description of acute cones of alcoves 
as intersections of half-spaces as in Lemma 5.8.
\end{proof}

As in section 2, define the extended affine Weyl group  
$\widetilde{W}^{[\Theta]} = X^{[\Theta]}_* \rtimes W^\Theta_0$, which acts on the set 
${\mathcal A}^{[\Theta]}$ of alcoves
of $V^{[\Theta]}$. Let $A'_0$ denote the base alcove and $\Omega^{[\Theta]}$ the stabilizer of 
$A'_0$ in $\widetilde{W}^{[\Theta]}$.  
Because $W^\Theta_{\rm aff}$ acts simply transitively on ${\mathcal A}^{[\Theta]}$, we have 
the decomposition $\widetilde{W}^{[\Theta]} = W^\Theta_{\rm aff} \rtimes \Omega^{[\Theta]}$.  
 
According to Proposition 9.3, we have $A_0 \cap V^{[\Theta]} = A'_0$, thus $\Omega\cap \widetilde{W}^{[\Theta]}
=\Omega^{[\Theta]}$. We also have $W_{\rm aff} \cap \widetilde{W}^{[\Theta]}=
W_{\rm aff}^{\Theta}$. This allows us to generalize the Kottwitz-Rapoport result on 
the inheritance property of the Bruhat order from affine Weyl groups to extended affine Weyl groups.

\begin{prop}
The Bruhat order $\leq$ on $\widetilde{W}^{[\Theta]}$ is inherited from the Bruhat order 
$\preceq$ on $\widetilde{W}$.  In other words, if $x,y \in \widetilde{W}^{[\Theta]}$, then 
$x \leq y$ if and only if $x \preceq y$.
\end{prop}

\begin{proof} According to Proposition 2.3 of Kottwitz-Rapoport \cite{KoRa}, the 
statement 
is already valid if $\widetilde{W}^{[\Theta]}$ is replaced by $W_{\rm aff}^{\Theta}$ and $\widetilde W$
is replaced by $W_{\rm aff}$. Since the Bruhat order on $W_{\rm aff}^{\Theta}$ (resp.
$W_{\rm aff}$) is extended in the obvious way to $\widetilde{W}^{[\Theta]}$ (resp. $\widetilde W$)
as explained in section 2, the proposition follows.
\end{proof}
 
Now let $\mu \in X^{[\Theta]}_*$ be a fixed coweight, which we suppose is dominant 
(note that it is dominant with respect to $R$ if and only if it is dominant with respect to 
$R^{[\Theta]}$).  We denote 
by ${\rm Adm}^\Theta(\mu)$, ${\rm Perm}^\Theta(\mu)$, and ${\rm Perm}^{st,\Theta}(\mu)$
the subsets of $\widetilde{W}^{[\Theta]}$ analogous to the subsets ${\rm Adm}(\mu)$, 
${\rm Perm}(\mu)$, and ${\rm Perm}^{st}(\mu)$ of $\widetilde{W}$.  The main goal of this 
section is the following proposition.  If $R$ is taken to be the root system for ${\rm GL}_{2n}$ and $\Theta$ its nontrivial automorphism (cf. section 10), then $R^{[\Theta]}$ is the root system for ${\rm GSp}_{2n}$.  Therefore this result implies Proposition 5.

\begin{prop}
Suppose that the root system $R$ is of type $A_m$. 
Then the equality
$${\rm Adm}^\Theta(\mu) = {\rm Perm}(\mu) \cap \widetilde{W}^{[\Theta]}$$
holds for every dominant coweight $\mu$ of $R^{[\Theta]}$.
\end{prop}

\begin{proof} Under the assumption that $R$ is of type $A_m$, we have
$${\rm Adm}(\mu)={\rm Perm}(\mu)={\rm Perm}^{st}(\mu).$$
If $x\in {\rm Adm}^\Theta(\mu)$, we have $x\in {\rm Adm}(\mu)$ according to Proposition 9.6,
thus $x\in {\rm Perm}(\mu) \cap \widetilde{W}^{[\Theta]}$.

Suppose $x \in {\rm Perm}(\mu) \cap \widetilde{W}^{[\Theta]}$.  Let $A = x(A_0)$ and 
$A' = x(A'_0)$; the alcoves $A$ in $V$ and $A'$ in $V^{[\Theta]}$ satisfy $A' = A \cap V^{[\Theta]}$. 
By Corollary 5.6, there is an element $w \in W^\Theta_0$ such that $A' \in C(A'_0,w)$.  
By Proposition 9.5, $A \in C(A_0,w)$.  
Since $x\in {\rm Perm}^{st}(\mu)$, we have $x(v) \in B(t_{w\mu}(v),w)$ for every vertex 
$v \in A_0$.  We can now apply Proposition 6.1. to $x$ and $t_{w\mu}$
and we obtain $x\leq t_{w\mu}$ relative to the Bruhat order on $\widetilde W$. 
By Proposition 9.6. the same inequality holds for the Bruhat order on $\widetilde{W}^{[\Theta]}$.
\end{proof}

Note that the set ${\rm Perm}(\mu) \cap \widetilde{W}^{[\Theta]}$ has some geometric 
meaning, when $m=2n-1$, which directly relates to the local models of \cite{H-N} attached to ${\rm GSp}_{2n}$ and a dominant coweight $\mu$.  More precisely, let us write $\Theta$ for the automorphism of ${\rm GL}_{2n}$ given by
$X \mapsto \tilde{J}^{-1} (X^t)^{-1}\tilde{J}$, where $\tilde{J}$ is the matrix
$$\begin{pmatrix} 0 & J \\ -J & 0 \end{pmatrix}, \quad$$
and where $J$ is the anti-diagonal matrix with entries equal to $1$.  This determines a symplectic group ${\rm Sp}_{2n} = {\rm SL}_{2n}^\Theta$, and a corresponding group ${\rm GSp}_{2n}$.
The automorphism $\Theta$ preserves the standard splitting for ${\rm GL}_{2n}$, hence induces an automorphism, also denoted $\Theta$, of the root system $R$ for ${\rm GL}_{2n}$.  The root system $R^{[\Theta]}$ of Steinberg is the root system of the group ${\rm GSp}_{2n}$ determined by its standard splitting (the ``upper triangular'' Borel subgroup and the diagonal torus).  

Fix a dominant cocharacter $\mu=(\mu_1,\dots,\mu_n,\nu_n,\dots,\nu_1)$ of ${\rm GSp}_{2n}$.  We can interpret the set ${\rm Perm}(\mu) \cap \widetilde{W}({\rm GSp}_{2n})$ in terms of lattice chains.  Let $k$ denote an algebraic closure of $\FF_p$, and let ${\mathcal V}_i = t^{-1}k[[t]]^i 
\oplus k[[t]]^{2n-i}$, for $0 \leq i \leq 2n$; extend by periodicity to get the ``standard'' infinite lattice chain ${\mathcal V}_\bullet$.  Consider the following set of periodic lattice chains
$$
\dots \subset {\mathcal L}_{-1} \subset {\mathcal L}_0 \subset \dots \subset {\mathcal L}_{2n} = t^{-1}{\mathcal L}_0 \subset \dots
$$
consisting of $k[[t]]$-lattices in $k(\!(t)\!)^{2n}$ with the following properties:

\begin{itemize}
\item ${\rm inv}({\mathcal L}_i,{\mathcal V}_i) \preceq \mu$, for every $i \in \ZZ$, 

\smallskip

\item ${\mathcal L}^{\perp}_i = t^{-c(\mu)}{\mathcal L}_{-i}$, for every $i \in \ZZ$.
\end{itemize}

Here ${\rm inv}(L,L')$ denotes the standard notion of invariant between two $k[[t]]$-lattices, and $\preceq$ denotes the usual partial order on dominant coweights for ${\rm GL}_{2n}$.  Moreover, $\perp$ is defined using the symplectic form $(x,y) \mapsto x^t\tilde{J}y$ on $k(\!(t)\!)^{2n}$, and $c(\mu)$ is the common value for $\mu_i + \nu_i$, $\,\, 1 \leq i \leq n$.   

This is precisely the set of $k$-points $M_{\mu}(k)$ for the local model $M_\mu$ considered in \cite{H-N}.  It carries an action of the standard Iwahori subgroup $I_k$ of 
${\rm GSp}_{2n}(k[[t]])$, namely the stabilizer of the standard lattice chain ${\mathcal V}_\bullet$.  

The $I_k$-orbits are clearly parametrized by the set ${\rm Perm}(\mu) \cap \widetilde{W}({\rm GSp}(2n))$.  The content of Proposition 9.7 is therefore that the strata in the special fiber of the model $M_\mu$ are indexed by the set of $\mu$-admissible elements of $\widetilde{W}({\rm GSp}_{2n})$.  As indicated by G\"{o}rtz \cite{Goertz1},\cite{Goertz2}, this makes it very likely that the models $M_\mu$ from \cite{H-N} are {\em flat}.

\section{Admissible and permissible sets for ${\rm GSp}_{2n}$}

In light of Theorem 3, it makes sense to ask when the equality ${\rm Adm}(\mu) = 
{\rm Perm}(\mu)$ holds for a dominant coweight $\mu$ of the root system for ${\rm GSp}_{2n}$.  
This root system arises as the ``fixed-points'' under an automorphism $\Theta$ of that of 
${\rm GL}_{2n}$, as in section 9.  By Proposition 9.7, we obviously have the following criterion:
$$
{\rm Adm}^\Theta(\mu) = {\rm Perm}^\Theta(\mu) ~ \Longleftrightarrow ~ {\rm Perm}^\Theta(\mu) \subset {\rm Perm}(\mu) \cap \widetilde{W}^{[\Theta]}.
$$
Using this will deduce that the equality ${\rm Adm}^\Theta(\mu) = {\rm Perm}^\Theta(\mu)$ 
holds whenever 
$\mu$ is a {\em sum of minuscule coweights} of ${\rm GSp}_{2n}$.

Define $X^* = X_* = \ZZ^{2n}$ and equip $V = \RR^{2n}$ with the standard inner product $(~,~)$.  
Let $R = \Check{R} = \lbrace e_i - e_j ~ | ~ 1 \leq i < j \leq 2n \rbrace$, where the $e_i$ are the 
standard basis vectors.  Let $\Pi = \lbrace e_i - e_{i+1} ~ | ~ 1 \leq i < 2n \rbrace$.  Then 
$R = (X^*,X_*,R,\Check{R},\Pi)$ is the root system for ${\rm GL}_{2n}$.  The finite Weyl group 
$W_0 =  S_{2n}$ acts on $X_* = \ZZ^{2n}$ by permuting the coordinates, and the extended 
affine Weyl group is $\widetilde{W} =  \ZZ^{2n} \rtimes S_{2n}$.  

Let $\Theta$  denote the automorphism of $X_* = \ZZ^{2n}$ defined by
$$
\Theta (x_1, \dots,x_{2n}) = (-x_{2n}, \dots,-x_1).
$$
We can then form the root system $R^{[\Theta]}$ as in section 9.  It is easy to check that 
$R^{[\Theta]}$ is the root system of ${\rm GSp}_{2n}$, so all the results of section 9 are in 
force in studying the latter.

\begin{theorem}
Let $\mu$ be a sum of minuscule coweights for ${\rm GSp}_{2n}$.  Then 
$$
{\rm Adm}(\mu) = {\rm Perm}(\mu).
$$
\end{theorem}

\noindent {\em Proof}.  According to Proposition  9.7,  we only need to show 
$$
{\rm Perm}^\Theta(\mu) \subset {\rm Perm}(\mu) \cap \widetilde{W}^{[\Theta]}.
$$  
The arguments we use are very close to those in Lemma 12.4 of \cite{KoRa}. 

We will need some notation.  For a vector $v \in \ZZ^{2n}$, let $v(m)$ denote its $m$th entry, so that $v = (v(1),\dots,v(2n))$.  Let $r = -\Theta$, that is, $r$ is the automorphism of $\RR^{2n}$ which reverses the order of the 
coordinate entries.  If $v,v' \in \RR^{2n}$, we write $v \leq v'$ if $v(m) \leq v'(m)$ for each $m$.  For $c \in \ZZ$, let ${\bf c} = (c^{2n}) = (c,\dots,c) \in \ZZ^{2n}$. 

Let $\omega_i = (1^i,0^{2n-i})$ for $1 \leq i \leq 2n$.  
The vectors $\omega_i$, $1 \leq i \leq 2n-1$, together with {\bf 0} serve as ``vertices'' of the base alcove $A_0$ for ${\rm GL}_{2n}$. 

The only minuscule coweights are the vectors $\nu$ in $\ZZ^{2n}$ of the form 
$\nu = (1^n,0^n) + {\bf c}$ or $\nu = {\bf c}$, for $c \in \ZZ$.  
Therefore we may write $\mu = (a^n,b^n)$, where $a,b \in \ZZ$ and $a \geq b$.  
For simplicity we discuss only the case where $b=0$, the more general case being similar. 

We choose a coordinate system $(x_1,\dots,x_n,y_n,\dots,y_1)$ on $V = \RR^{2n}$.  Then 
$X^{[\Theta]}_* = \lbrace (x_1,\dots,x_n,y_n,\dots,y_1) \in \ZZ^{2n} ~ | ~ x_1+y_1 = \cdots =x_n + y_n \rbrace$.
The finite Weyl group $W^{[\Theta]}_0$ is $W_n = (\ZZ/2\ZZ)^n \rtimes S_n$.  
The element $e_i \in (\ZZ/2\ZZ)^n$ acts by switching $x_i$ and $y_i$, and an element $\sigma \in S_n$ acts by simultaneously permuting the $x_i$'s and the $y_i$'s.
Let ${\rm Conv}^\Theta(\mu)$ denote the convex hull in $\RR^{2n}$ of the set $W_n(\mu)$.  We have 
$$\begin{array}{lll}
{\rm Conv}(\mu) &=& \lbrace (x_1,\dots,x_n,y_n,\dots,y_1) \in \RR^{2n} ~ | ~ 0 \leq x_i,y_i  \leq a, \forall i, \,\, \mbox{and} \,\, \sum_i x_i + y_i = na \rbrace \\
&=&  \lbrace v \in \RR^{2n} ~ | ~ {\bf 0} \leq v \leq {\bf a} \,\, \mbox{and} \,\, \sum_m v(m) = na \rbrace,
\end{array}$$
and
$$\begin{array}{lll}
{\rm Conv}^\Theta(\mu)&=& \lbrace (x_1,\dots,x_n,y_n,\dots,y_1) \in \RR^{2n} ~ | ~ 0 \leq x_i,y_i \leq a \,\, \mbox{and} \,\, x_i + y_i = a, \forall i \rbrace \\
&=& \lbrace v \in \RR^{2n} ~ | ~ {\bf 0} \leq v \leq {\bf a} \,\, \mbox{and} \,\, v + r(v) = {\bf a} \rbrace.
\end{array}$$

The vectors $\eta_i := (\omega_i + \omega_{2n-i})/2$, $1 \leq i \leq n$, together with {\bf 0}, serve as ``vertices'' for the base alcove $A'_0$.  
We now fix $x \in \widetilde{W}^{[\Theta]}$, and let $v_i := x(\omega_i)$, for $1 \leq i \leq 2n-1$.  One can easily verify that $x(\eta_i) - \eta_i = (v_i - \omega_i + v_{2n-i} - \omega_{2n-i})/2$, for $1 \leq i \leq n$.

Now suppose $x = t_\lambda w \in {\rm Perm}^\Theta(\mu)$, i.e., 
$\lambda = x({\bf 0})-{\bf 0}$ and $x(\eta_i)-\eta_i$ belong to ${\rm Conv}^\Theta(\mu)$ for $1 \leq i \leq n$.  Therefore
\begin{gather}
{\bf 0}  \leq \lambda \leq {\bf a} \tag{10.1} \\
\lambda + r(\lambda)   =  {\bf a} \tag{10.2} \\ 
{\bf 0}  \leq  v_i - \omega_i + v_{2n-i} - \omega_{2n-i}  \leq  {\bf 2a} \tag{10.3} 
\\
v_i - \omega_i + v_{2n-i} - \omega_{2n-i} + r(v_i - \omega_i + v_{2n-i} - 
\omega_{2n-i}) = {\bf 2a}, \tag{10.4} 
\end{gather}
the last two equations holding for all $i=1,\dots,n$.  We also have

\begin{gather}
v_i + r(v_{2n-i}) = {\bf a+1} \tag{10.5} \\
\omega_i + r(\omega_{2n-i}) = {\bf 1}, \tag{10.6}
\end{gather} 
for all $i=1,\dots, 2n$  (for (10.5), use (10.2) and (10.6) together with the observation that 
$w \in W^\Theta_0 \Rightarrow w \,\, \mbox{commutes with} \,\, r$). 

We need to show that $x \in {\rm Perm}(\mu)$, i.e., $\lambda$ and $x(\omega_i)-\omega_i$ belong to ${\rm Conv}(\mu)$ for $1 \leq i \leq 2n-1$.  Thus we need to show

\begin{gather}
{\bf 0} \leq \lambda \leq {\bf a} \tag{10.7} \\
\sum_m \lambda(m) = na \tag{10.8} \\
{\bf 0} \leq v_i - \omega_i \leq {\bf a} \tag{10.9} \\
\sum_m v_i(m) - \omega_i(m) = na \tag{10.10},
\end{gather}
the last two equations holding for all $i = 1,\dots,2n-1$.  
It is easy to check that (10.7), (10.8), and (10.10) are satisfied.  
It is enough to verify (10.9) for $i$ such that $1 \leq i \leq n$: applying $r$ 
to (10.9) and using the identities in (10.5-10.6) yields (10.9) with $2n-i$ 
replacing $i$.  Henceforth we fix $i$ in this range.  From (10.5-10.6) we can 
easily verify that ${\bf a} - r(v_i - \omega_i) = v_{2n-i} -\omega_{2n-i}$.  
Therefore (10.2) and (10.3), respectively, yield 

\begin{gather}
{\bf a} \leq v_i + r(v_i) \leq {\bf a+1} \tag{10.11} \\
{\bf -a} \leq v_i-\omega_i - r(v_i-\omega_i) \leq {\bf a}. \tag{10.12}
\end{gather}

\noindent Let $u = v_i(m) - \omega_i(m)$ and $v = r(v_i)(m) - r(\omega_i)(m)$ for any $m$ in the range $1 \leq m \leq 2n$.  Then (10.11) and (10.12) yield
\begin{gather}
a \leq u + v + \delta \leq a+1 \tag{10.13} \\
 -a \leq u-v \leq a, \tag{10.14} 
\end{gather}
where $\delta = \omega_i(m) + r(\omega_i)(m)$.  Since $1 \leq i \leq n$, the term $\delta$ is either 0 or 1.  Adding (10.13) and (10.14) gives
\begin{equation}
0 \leq 2u + \delta \leq 2a +1 \tag{10.15}
\end{equation}
and thus $0 \leq u = v_i(m) - \omega_i(m) \leq a$.  Since this inequality holds for every $m$, we have ${\bf 0} \leq v_i - \omega_i \leq {\bf a}$, as desired.

\hfill $\square$

\section{Remarks on ${\rm PSO}(2n+1)$}

The results we obtained in the symplectic case depend in a crucial way on the fact that
the restriction of the Bruhat order of $\widetilde W({\rm GL}(2n))$ to $\widetilde W({\rm GSp}(2n))$
is the Bruhat order of the latter group. This inheritance property is no longer true for the odd 
orthogonal group. We are grateful to Kottwitz for his help in preparing this section. 

The group ${\rm PSO}(2n+1)$ whose root system is $B_n$ 
can be defined as the fixed point subgroup of an involution of 
${\rm PGL}(2n+1)$. This induces an involution $\Theta$ on the root system $A_{2n}$. 
However the fixed point root system of $\Theta$ in the sense of Steinberg, see \cite{Steinberg} 
or section 9, is not $B_n$ but rather $C_n$. We have the natural inclusions
$$W_{\rm aff}(B_n)\subset W_{\rm aff}(C_n)\subset W_{\rm aff}(A_{2n}).$$
According to \cite{Steinberg} and \cite{KoRa}, the Bruhat order of $W_{\rm aff}(C_n)$ is inherited
from $W_{\rm aff}(A_{2n})$. We will show by the following examples that the Bruhat order of
$W_{\rm aff}(B_n)$ is not.

Let denote $s_0,\ldots,s_n$ the simple reflections of $W_{\rm aff}(B_n)$ and 
$s'_0,\ldots,s'_n$ those of $W_{\rm aff}(C_n)$. The inclusion $W_{\rm aff}(B_n)\subset 
W_{\rm aff}(C_n)$ is given by $s_i \mapsto s'_i$ for $i\not=0$ and $s_0\mapsto s'_0 s'_1s'_0$.
The elements $s_0$
and $s_1$ are not related by the Bruhat order in $W_{\rm aff}(B_n)$ but their images 
$s'_0s'_1s'_0$ and $s'_1$ are related. 

\bigskip

\begin{center} 
\includegraphics{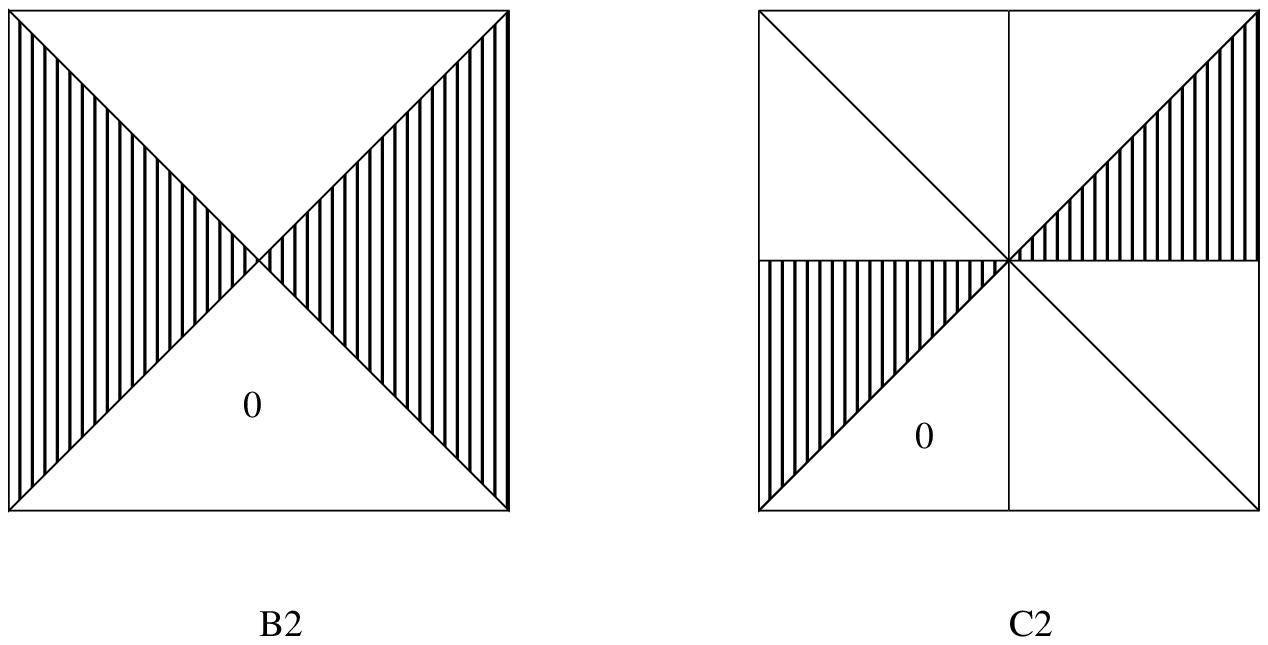}
\end{center}

\bigskip
The extended affine Weyl group $\widetilde W(B_n)=W_{\rm aff}(B_n) \rtimes \{1,\tau\}$ 
is canonically identified with $W_{\rm aff}(C_n)$ by $\tau\mapsto s'_0$. One can show 
by similar examples that the extended Bruhat order on the coset $W_{\rm aff}(B_n)\tau$
is also not inherited from $W_{\rm aff}(C_n)$.

Let $\mu \in \ZZ^n$ be a dominant coweight of type $B_n$.  Concerning the analog of Theorem 4 
for ${\rm PSO}(2n+1)$, we do not know whether the equality
$$
{\rm Adm}^{B_n}(\mu) = {\rm Perm}^{B_n}(\mu)$$
holds for $\mu$ a sum of minuscule coweights.  However, it is not hard to see that the analog 
of Proposition 5 is false.  

\smallskip

\noindent {\bf Example.}
Let $n=2$ and $\mu = (1,0)$.  Then ${\rm Adm}^{B_n}(\mu)$ consists of 13 elements.  
But 
$${\rm Perm}^{A_{2n}}(\mu) \cap \widetilde{W}(B_n) = {\rm Perm}^{A_{2n}}(\mu) \cap W_{\rm aff}(C_n) 
= {\rm Adm}^{C_n}(\mu)$$
consists of 19 elements.  The final equality above follows from Proposition 9.7. 

\medskip

\section{Acknowledgments}
We thank R. Kottwitz for encouragement, helpful comments, and a careful reading of the paper.  
We are also grateful to M. Rapoport for communicating his guess that Theorem 4 might hold, and 
for his continued interest in this work.  Finally, we thank W. Casselman for his comments on the paper.

This article was written during the {\em Semestre Hecke} workshop on automorphic forms at the 
Institut Henri Poincar\'{e} in Paris, held in the spring semester of 2000.  The first author thanks the 
IHP for its hospitality.  His research is partially supported by an NSF-NATO postdoctoral fellowship, 
an NSERC research grant, a Connaught New Staff grant, and NSF grant DMS 97-29992. The 
second author would like to thank the mathematics department of the University of Toronto 
for a fruitful visit in February 2000 during which this research in collaboration could begin.

\small
\bigskip
\obeylines
\noindent
Institute for Advanced Study
School of Math
Einstein Drive
Princeton, NJ 08540
email: haines@ias.edu

\bigskip
\obeylines
CNRS, Universit\'{e} de Paris Nord
D\'epartement de math\'ematiques, 
93430 Villetaneuse, FRANCE
email: ngo@math.univ-paris13.fr


\begin{thebibliography}{20s}


\bibitem{BGG1} I.N. Bernstein, I.M. Gelfand, and S.I. Gelfand, {\it Differential operators on the base affine space and a study of ${\mathfrak g}$-modules}, Prepr. 77, IPM Akad. Nauk SSSR (1972) (English translation in: Lie groups and their representations. Proc. Summer School in Group Representations. Bolyai Janos Math. Soc., Budapest 1971, pp.21-64, New York: Halsted (1975)).

\bibitem{BGG} I.N. Bernstein, I.M. Gelfand and S.I. Gelfand, 
{\it Schubert cells and cohomology of the spaces $G/P$}, Uspekhi Mat. Nauk. 
{\bf 28} (1973), 3-26; English translation, Russ. Math. Surv. {\bf 28} (1973), 1-26.

%\bibitem{Deod}
%V. Deodhar, {\em On some geometric aspects of Bruhat 
%orderings. I. A finer decomposition of Bruhat cells}, 
%Invent. Math. {\bf 79}, 499-511 (1985).

\bibitem{Deod1}
V. Deodhar, {\em Some characterizations of the Bruhat ordering on a Coxeter 
group and determination of the relative M\"{o}bius function}, 
Invent. Math. {\bf 39} (1977), 187-198.

\bibitem{Deod2}
V. Deodhar, {\em On Bruhat ordering and weight-lattice ordering for a Weyl group}, 
Indag. Math. {\bf 40} (1978), 423-435.

%\bibitem{Del-Rap}
%P. Deligne and M. Rapoport, {\it Modules des courbes elliptiques}, in: Modular functions of one variable, II, Springer Lecture Notes {\bf 349} (1973).


\bibitem{Goertz1}
U. G\"{o}rtz, {\em On the flatness of certain Shimura varieties of PEL-type}, preprint 
K\"{o}ln 1999.

\bibitem{Goertz2}
U. G\"{o}rtz, forthcoming paper.

\bibitem{Haines}
T. Haines, {\em The Combinatorics of Bernstein functions}, preprint (1998), {\em to appear in} Trans. Amer. Math. Soc.

\bibitem{H-N}
T. Haines and B.C. Ng\^{o}, {\em Nearby cycles for local models of some Shimura varieties}, preprint (1999).

\bibitem{Hum}
J.E. Humphreys, {\it Reflection Groups and Coxeter Groups}, Cambridge Studies in Advanced Mathematics, no. 29, Cambridge Univ. Press (1990).

\bibitem{IM}
N. Iwahori and H. Matsumoto, {\em On Some Bruhat Decomposition 
and the Structure of the Hecke Rings of $p$-adic Chevalley Groups}, 
Publ. Math. IHES, No. 25, (1965) 5-48.


%\bibitem{Lusz}
%G. Lusztig, {\em Singularities, character formulas and a $q$-analog of weight multiplicities}, Ast\'{e}risque {\bf 101-102} (1983), 208-229.

%\bibitem{AffineHecke}
%G. Lusztig, {\em Affine Hecke Algebras and Their Graded Version}, 
%J. of the Amer. Math. Soc. {\bf 2} No. 3 (1989) 599-635.

\bibitem{KoRa}
R. Kottwitz, M. Rapoport,  {\em Minuscule  Alcoves for $Gl_n$ and $GSp_{2n}$},
manuscripta math. {\bf 102}, 403-428 (2000). 

\bibitem{RapAA}
M. Rapoport, {\em On the Bad Reduction of Shimura Varieties}, 
Automorphic Forms, Shimura Varieties and
$L$-functions, part II, Perspectives in Mathematics, vol. 11, Academic
Press, San Diego, CA, 1990, 253-321. 

\bibitem{Rap-Zink} 
M.~Rapoport,Th.~ Zink.\newblock { Period spaces for $p$-divisible groups}.
\newblock {\em Annals of Math. Studies} 144, Princeton Univ. Press 1996

\bibitem{Springer}
T.A. Springer, {\em Reductive groups}, Automorphic forms, Representations, and $L$-functions, part I, Proc. of the Symp. in Pure Math., vol. XXXIII, 
Amer. Math. Soc. 1977.

\bibitem{Steinberg}
R. Steinberg, {\em Endomorphisms of linear algebraic groups}, Mem. Amer. Math. Soc. {\bf 80} (1968), 1-108.

%\bibitem{Rapletter}
%M. Rapoport, {\em Letter to Waldspurger}, January/February 1989.

%\bibitem{Waldspurger}
%J.L. Waldspurger {\em Letter to Rapoport}, February 1989.


\bigskip

\end{thebibliography}
\end{document}